\documentclass[11pt]{article}   
\usepackage{tikz}
\usepackage{tikz-3dplot}
\usetikzlibrary{hobby}
\usetikzlibrary{positioning}
\usetikzlibrary{shapes}
\usepackage{graphicx}
\usepackage{color}
\usepackage{subfigure}
\usepackage{array}
\usepackage{longtable}

\usepackage{amsmath, amssymb, amsthm}
\usepackage{graphicx}
\usepackage{hyperref}
\usepackage{geometry}
\usepackage{cite}

\usepackage{authblk} 
\usepackage{abstract} 
\usepackage{enumitem} 

\title{Hybrid PDE-Deep Neural Network Model for Calcium Dynamics in Neurons}

\author[ ]{Abel Gurung 
	\thanks{This work was completed while the first author was an undergraduate student at the University of Southern Mississippi. Email: \texttt{abel.gurung@usm.edu}}}
\author[ ]{Qingguang Guan \thanks{Correspondence to: \texttt{qingguang.guan@usm.edu}}}
\affil[ ]{School of Mathematics and Natural Sciences\\ 
	University of Southern Mississippi\\
	 118 College Drive, Hattiesburg, MS, 39406}
	 
\date{}

\begin{document}
	
	\maketitle

\begin{abstract}
Traditionally, calcium dynamics in neurons are modeled using partial differential equations (PDEs) and ordinary differential equations (ODEs). The PDE component focuses on reaction-diffusion processes, while the ODE component addresses transmission via ion channels on the cell's or organelle's membrane. However, analytically determining the underlying equations for ion channels is highly challenging due to the complexity and unknown factors inherent in biological processes. Therefore, we employ deep neural networks (DNNs) to model the open probability of ion channels, a task that can be intricate when approached with ODEs. This technique also reduces the number of unknowns required to model the open probability. When trained with valid data, the same neural network architecture can be used for different ion channels, such as sodium, potassium, and calcium. Furthermore, based on the given data, we can build more physiologically reasonable DNN models that can be customized. Subsequently, we integrated the DNN model into calcium dynamics in neurons with endoplasmic reticulum, resulting in a hybrid model that combines PDEs and DNNs. Numerical results are provided to demonstrate the flexibility and advantages of the PDE-DNN model.
\end{abstract}


\section{Introduction}\label{sec1}
Calcium signals (changes in the concentration of Ca$^{2+}$ ions) play essential roles in cells, such as 
neurons, muscle and immune cells. In neurons, calcium signals regulate multiple processes including neurotransmitter release (see \cite{augustine2001does,sudhof2012calcium}) and neuronal health (see \cite{ermak2002calcium, chen2005ca2+,jia2023change}). When an action potential (the electrical signal that involves a change in voltage across the cell membrane) reaches the presynaptic terminal of a neuron, voltage-gated calcium channels open, allowing Ca$^{2+}$ ions to flow into the active zone. The influx of Ca$^{2+}$ triggers the fusion of synaptic vesicles with the presynaptic membrane, resulting in the release of neurotransmitters into the synaptic cleft. Calcium signaling also plays a dual role in neuroprotection and neurotoxicity. Moderate and regulated calcium signals are involved in protective pathways that maintain neuronal health. However, excessive calcium influx, such as during excitotoxic events (e.g., stroke or trauma), can lead to cell damage and death through the activation of destructive enzymes like calpains. Calcium signals are just as crucial for muscle cells as they are for neurons. These signals are essential for both initiating and terminating muscle contraction (see \cite{hill2011calcium,berridge2008smooth,wray2005calcium}). In response to an action potential, Ca$^{2+}$ release from the sarcoplasmic reticulum (SR - an organelle in muscle cells used to store calcium ions) triggers muscle contraction. Muscle relaxation occurs when Ca$^{2+}$ is pumped back into the SR. This precise regulation ensures proper muscle function during activities such as movement and heartbeat. In T cell activation (see \cite{oh2008calcium,vig2009calcium,trebak2019calcium}), antigen recognition by the T cell receptor triggers intracellular signaling pathways, leading to calcium release from intracellular stores (such as ER - endoplasmic reticulum). Calcium ions activate proteins like calcineurin.  Calcineurin, in turn, influences gene expression, which regulates T cell activation and effector functions. Calcium signaling thus plays a crucial role in coordinating T cell responses to antigens presented by antigen-presenting cells. 

\begin{figure}[ht!]   
	\begin{center}
		\includegraphics[width=0.6\textwidth]{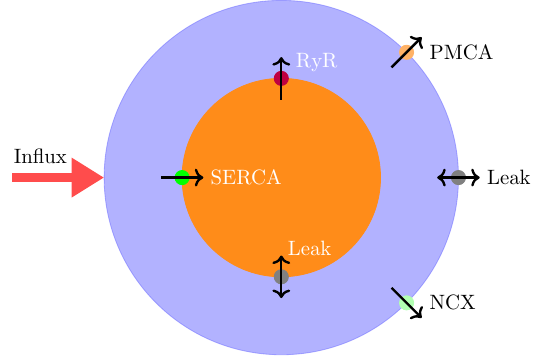}
	\end{center}
	\caption{Cross-section of the axon (It can also be viewed as a two-dimensional closed cell). The orange region is ER, denoted as $\Omega_e$. The blue region is cytosol, denoted as $\Omega_c$. There are no other organelles inside the axon. ER membrane is denoted as $\Gamma$, plasma membrane is denoted as $\partial\Omega$, where $\Omega = \Omega_e\cup\Gamma\cup \Omega_c$.}
	\label{fig1}
\end{figure}

Calcium channels on the membranes of cells and organelles are key to responding to stimuli like action potentials and generating calcium signals. There are many different types of calcium channels, such as
Voltage-Gated Calcium Channels (\cite{catterall2011voltage}), Ligand-Gated Calcium Channels (\cite{burnashev1998calcium,hucho2001ligand}) and Store-Operated Calcium Channels (\cite{prakriya2015store,croisier2013activation,luik2008oligomerization}). Those channels allow calcium flux through membranes, thus changing the calcium concentration in the cytosol or organelles. We focus on Ligand-Gated Calcium Channels on neurons, such as PMCA pump, NCX pump, and leak channel on the plasma membrane, and SERCA pumb, RyR channel, and leak channel on the ER membrane, see Figure \ref{fig1} for a 2D closed cell example. Calcium pumps move Ca$^{2+}$ ions from regions of low concentration (the cytosol) to regions of high concentration (extracellular region and ER). PMCA stands for Plasma Membrane Calcium ATPase, see \cite{carafoli1991calcium}, it is a type of enzyme found in the cell membrane that actively transports calcium ions out of the cell. The calcium flux through PMCA pump can be modeled by a second order hill function as in equation \eqref{pmca-ncx}, (see \cite{graupner2005theory}). NCX pump is the Na$^{+}$/Ca$^{2+}$ exchanger, which pumps out calcium ions in exchange for sodium ions coming into the cell. However, we will only consider the calcium ions. Thus, the flux pumped by NCX can be described by the first-order Hill function as in equation \eqref{pmca-ncx}, see \cite{breit2018}. The Sarco/Endoplasmic Reticulum Ca$^{2+}$-ATPase (SERCA) pump is an ATP-dependent enzyme that transports calcium ions from the cytosol into the sarcoplasmic or endoplasmic reticulum (SR/ER), see \cite{lipskaia2009role}. The mathematical model of the Ca$^{2+}$ flux crossing the SERCA pump in equation \eqref{ryr-serca}  is from \cite{sneyd2003model}, it depends on both cytosol and ER calcium concentrations. Leak channels allow exchange of Ca$^{2+}$ across cell and organelle membranes (see \cite{samtleben2015store,schulte2022shaped}). They are essential components of cellular membranes, contributing to reaching equilibrium across the membrane for calcium ions. The calcium flux is usually assumed to be proportional to the difference of calcium concentrations on both sides of the membrane.
All previously mentioned calcium channels can be modeled as functions dependent on the calcium concentration, independent of the dynamics (i.e., behavior of the concentration over time) of the calcium concentration. The RyR channel ( Ryanodine Receptor channel) is special because it depends on the dynamics of calcium concentration. We consider this channel separately.

RyR channels are embedded in the membranes of the ER (see \cite{santulli2018ryanodine,des2016structural,bezprozvanny1991bell}). They play a crucial role in regulating the release of calcium ions from the ER into the cytosol. The process is dynamic because a small rise in cytosolic calcium can promote further release of calcium from the ER through RyR channels.
Though it's not completely understood, RyR channels are assumed to be regulated solely by calcium in the cytosol from previous research, i.e., through calcium-induced calcium release. RyR channels are measured by their open probability, which refers to the likelihood that an RyR channel is open at any given moment. There are a lot of uncertainties and many different models about RyR channels.
We choose the model of open probability from \cite{keizer1996}, where the open probability $P$ ranges from 0 to 1. It depends on the time $t$ and calcium concentration $u$ in the cytosol:
\begin{equation}\label{ryr_ptu} 
	{P(t,u)} = 1-c_1(t,u)-c_2(t,u),
\end{equation}
$c_1(t,u), c_2(t,u)$ are functions of $t,u$, and calculated from an ODE system:
\begin{equation} \label{ryr_ode}
	{
		\begin{bmatrix}
			c_1'\\
			o'\\
			c_2'
		\end{bmatrix}
		=
		\begin{bmatrix}
			-u^4k_a^+-k_a^- & -k_a^- & -k_a^- \\
			-u^3k_b^+ & -u^3k_b^+-k_b^- & -u^3k_b^+ \\
			-k_c^+    & -k_c^+       & -k_c^+-k_c^- 
		\end{bmatrix}
		\begin{bmatrix}
			c_1\\
			o\\
			c_2
		\end{bmatrix}
		+
		\begin{bmatrix}
			k_a^-\\
			u^3k_b^+\\
			k_c^+
		\end{bmatrix}
	}
\end{equation}
where $k_a^{\pm},k_b^{\pm},k_c^{\pm}$ are positive constants obtained from fitting the experimental data. This model is derived from a four-state continuous-time Markov chain in \cite{keizer1996}, see \cite{Norris_1997} for an introduction to Markov chains. To get $P$, three extra unknowns ($c_1,o,$ and $c_2$) need to be solved.

Axons are typically long and slender, resembling a thread or wire in neurons. Unlike dendrites, which often branch, axons maintain a relatively constant diameter for most of their length. The ER in axons primarily exists as a tubular network, with these thin tubes extending throughout the length of the axon (see \cite{ozturk2020axonal} for an illustration).  We consider a simplified calcium model on axons from \cite{breit2018}.  

Our major contribution is proposing a deep neural network model to replace the existing one \eqref{ryr_ptu}-\eqref{ryr_ode} and building a better, more flexible hybrid PDE-DNN model.
 The DNN model can reduce the three unknowns of the ODE system \eqref{ryr_ode} to one. It is also more flexible and can be customized using observational or artificial data. Details of the drawbacks of the ODE model will be provided in the following sections. 
 
The paper is organized as follows:
In Section 2, we transform the PDE part of the calcium model to cylindrical coordinates and further reduce it to a one-dimensional reaction-diffusion system. In Section 3, we describe the structure of the DNN model and train it using data generated from the ODE model. We demonstrate that the DNN trained on this dataset inherits drawbacks from the ODE model. To address these issues, we propose training the DNN on reasonable artificial datasets, leading to a more flexible model for the RyR channel. Section 4 is devoted to the numerical methods for solving the PDE-DNN system. In Section 5, we present numerical results and discuss the disadvantages of the PDE-ODE model and the advantages of the PDE-DNN model when the DNN is trained on specially designed artificial datasets. Conclusions are drawn in Section 6.

\section{Calcium Model in Cylindrical Coordinates }
The axon is assumed to be a cylinder within a cylinder, both sharing the same axis. A cross-section of the axon is shown in Figure \ref{fig1}. Suppose that $u$ is the calcium concentration in cytosol ($\Omega_c$) and $u_e$ is the calcium concentration in ER ($\Omega_e$). If $b$ is the concentration of buffered calcium ions in $\Omega_c$ and the diffusivities of $u, u_e$ and $b$ are constants in $\Omega_c$ and $\Omega_e$, then we have a diffusion reaction system, which are given in equation \eqref{diff}
\begin{equation}\label{diff}
	\left\{
	\begin{aligned}
		&{\partial_t u} -\nabla\cdot ({D_c} \nabla u)  =  K_b^-(b^0-b)-K_b^+ bu
		      \quad \text{on } \Omega_c \qquad\qquad\\
  	&{\partial_t b} -\nabla\cdot ({D_b} \nabla b)  = K_b^-(b^0-b)-K_b^+ bu
				\quad \text{on } \Omega_c\\
		& {\partial_t u_e} -\nabla\cdot ({D_c} \nabla u_e)  = 0  \quad \text{on } \Omega_e
	\end{aligned}
	\right. 
\end{equation}
where $D_c, D_b$ are diffusion coefficients, $K_b^{+}, K_b^{-}$ are reaction coefficients, $b_0$ is the highest concentration  of unbinding calcium buffer. The boundary conditions are of Neumann type:
\begin{equation}\label{neu}
	\left\{
\begin{aligned}
	&D_b\partial_n b 
	= 0
	\ {\text{ on } \partial\Omega_c} ,\\
	&D_c\partial_n u 
	= g_c(u)
	\ {\text{ on } \partial\Omega} ,\\
	&D_c\partial_n u  
	= -g_e(u,u_e)
	\ \text{ on } \Gamma ,\\
	&D_c\partial_n u_e
	= g_e(u,u_e)
	\ \text{ on }  \Gamma,
	\qquad\qquad\qquad\qquad\qquad\qquad
\end{aligned}
	\right. 
\end{equation}
where $\partial \Omega_c = \Gamma\cup \partial\Omega$ and the calcium fluxes through plasma/ER membrane are
\begin{align}
	&g_c(u) =  
	-
	\underbrace{C^c_1\frac{u^2}{K_p^2+u^2}}_\text{PMCA}
	-
	\underbrace{C^c_2\frac{u}{K_n+ u}}_\text{NCX}
	+
	\underbrace{C^c_3(c_o-u)}_\text{Leak} \label{pmca-ncx}\\
	&g_e(u,u_e) =  
	\underbrace{C_{1}^e { P(t,u)} (u-u_e)}_\text{RyR}
	+
	\underbrace{C_2^e\frac{u}{(K_s+u)u_e}}_\text{SERCA}
	-
	\underbrace{C_3^e(u_e-u)}_\text{Leak} \label{ryr-serca}     
\end{align}
the coefficients $C_1^c,C_2^c,C_3^c$,$K_n, K_p, K_s$, $C_1^e,C_2^e,C_3^e$ are positive constants. Equations \eqref{ryr_ptu} to \eqref{ryr-serca} complete the model in Cartesian coordinates. The initial values for $u, b, u_e$, $c_1,o,c_2$ are all constants and make fluxes zero.
The well-posedness and efficient numerical methods for solving the model can be found in \cite{guan2022}.

Suppose that in a 3D Cartesian coordinate system, the axon's axis coincides with the $z$-axis. The radius of the endoplasmic reticulum (ER) is denoted by $L$, and the radius of the axon is denoted by $R$. Then, we can transform equation \eqref{diff} into cylindrical coordinates $(r, \theta, z)$.
To simplify, let's assume $u$, $b$, and $u_e$ have no dependence on $\theta$ and $z$. Under this assumption, equation \eqref{diff} becomes:
\begin{equation}\label{diff_r}
\left\{
\begin{aligned}
&\frac{\partial u}{\partial{t}} -D_c\frac{\partial^2 u}{\partial{r}^2}
-D_c\frac{1}{r}\frac{\partial u}{\partial{r}}  
= K_b^-(b^0-b)-K_b^+ bu 
\quad \text{on } (L,R), \qquad\\
&\frac{\partial b}{\partial{t}} -D_b\frac{\partial^2 b}{\partial{r}^2}
-D_b\frac{1}{r}\frac{\partial b}{\partial{r}}  
= K_b^-(b^0-b)-K_b^+ bu 
\quad \text{on } (L,R), \qquad\\
&\frac{\partial u_e}{\partial{t}} -D_c\frac{\partial^2 u_e}{\partial{r}^2}
-D_c\frac{1}{r}\frac{\partial u_e}{\partial{r}}  
= 0 
 \quad \text{on } (0,L).
\end{aligned}
\right. 
\end{equation}
The boundary conditions in \eqref{neu} become:
\begin{equation}\label{neu_r}
	\left\{
\begin{aligned}
	&D_b\partial_r b 
	= 0,
	\ {\text{ if }} r=R {\text{ or }}  r=L, \\
	&D_c\partial_r u 
	= g_c(u),
	\ {\text{ if }} r=R,\\
	&-D_c\partial_r u  
	= -g_e(u,u_e),
	\ \text{ if } r=L ,\\
	&D_c\partial_r u_e
	= g_e(u,u_e),
	\ \text{ if }  r=L,
	\qquad\qquad\qquad\qquad\qquad\qquad
\end{aligned}
	\right. 
\end{equation}
where the RyR model in equations \eqref{ryr_ptu}-\eqref{ryr_ode}, $g_c(u)$ and $g_c(u,u_e)$ in equations \eqref{pmca-ncx}-\eqref{ryr-serca} keep the same. However, we have $\partial_nu=-\partial_r u$ at $r=L$.  After simplification of the calcium model, we only need to deal with a 1D diffusion reaction system and solve the RyR model at a single radial position, $r = L$, rather than for all points on the ER membrane. 

For the simulations in the following sections, we will use a combined set of equations. This set includes:
\begin{itemize}
    \item  Equations \eqref{diff_r}--\eqref{neu_r} for the diffusion reaction process of calcium ions in cytosol and ER.
    \item Equations \eqref{ryr_ptu}--\eqref{ryr_ode} for the open probability of a single RyR channel. A novel substitution, modeled by a DNN, is also proposed in this work.
    \item Equations \eqref{pmca-ncx}--\eqref{ryr-serca} for calcium fluxes across the plasma and ER membranes.
\end{itemize}

\section{DNN Model for the Open Probability of RyR Channel}\label{sec3}
In the model \eqref{ryr_ptu}-\eqref{ryr_ode} for the RyR channel, the open probability $P(t)$ depends on $u(t)$ (the calcium concentration), and two intermediate values, $c_1(t)$ and $c_2(t)$. However, $c_1(t)$ and $c_2(t)$ are solved from the ODE system with three unknowns. Since we only need $P(t)$ in the calcium model, is it possible to model $P(t)$  directly without the intermediate values? Of course, this would require a new model, which must be significantly different but still capable of capturing solutions similar to the old one. Here we propose a possible model with $P(t)$ as the only unknown:
\begin{equation}
	\left\{
\begin{aligned}
    &\frac{{\rm d} P}{{\rm d} t} = F(P(t),u(t),u'(t))\\
    &P(0)=p_0
\end{aligned} 
\right.
\end{equation}
where $u(t)$ is a given function of the calcium concentration, which varies as time changes, $u'(t)$ is the time derivative of $u(t)$, $p_0$ is the initial open probability, but $F$ is an unknown scalar-valued function with three variables. Now the problem becomes: how do we determine $F$?

Fully connected deep neural networks have been successfully applied to model ODEs; see \cite{su2021deep, qin2021data, churchill2023flow, chen2021generalized} for an introduction to flow-map learning for dynamical systems. In \cite{lei2021neural}, flow-map learning is applied to model an ion channel with one unknown. However, we are not learning the flow-map of the ODE system in equation \eqref{ryr_ode}. Instead, we are reducing the order of the ODE system from three to one. This strategy can also be applied to reduce the order of other continuous-time Markov chain models for ion channels with more unknowns, such as \cite{wang1997quantitative,whittaker2022ion}.  

\begin{figure}[ht!]  
	\begin{center}
		\includegraphics[width=0.9\textwidth]{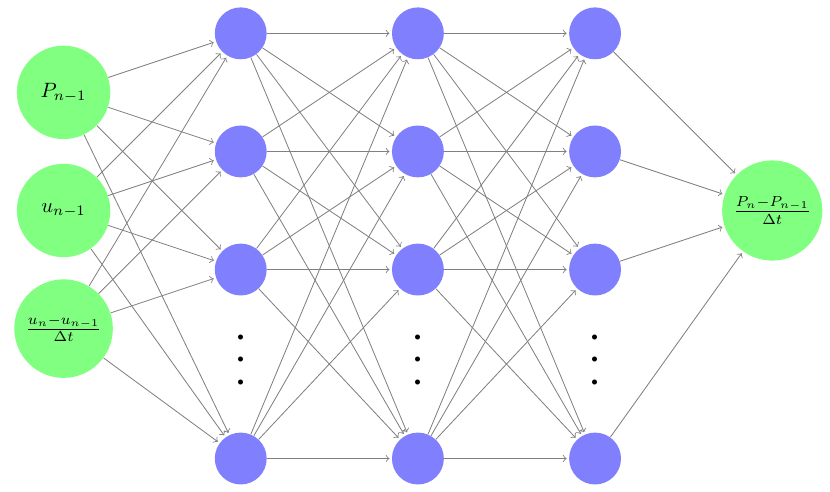}
	\end{center}
	\caption{The structure of the fully connected DNN model for the RyR channel. There are three inputs, one output, and three hidden layers. Among the inputs, $P_{n-1}$ is the open probability of the RyR channel at time $t_{n-1}$, $u_{n-1}$ is the calcium concentration in the cytosol at $t_{n-1}$, and similarly, $u_n$ is the Ca$^{2+}$ concentration at $t_n$. $(u_n - u_{n-1})/\Delta t$ approximates the time derivative of $u$ at $t_{n-1}$. The output $(P_n - P_{n-1})/\Delta t$ is an approximation of the time derivative of $P$ at $t_{n-1}$.}
	\label{fig2}
\end{figure}
\subsection{The Structure of the Fully Connected DNN Model}
We propose a fully connected deep neural network to determine $F$, using a structure that consists of three hidden layers, see Figure \ref{fig2}. The inputs are designed as follows:
\begin{enumerate}
	\item \textbf{Probability of Ion Channel Opening:} The first input $P_{n-1}=P(t_{n-1})$ is the open probability of RyR channel at time $t_{n-1}$. 
	\item \textbf{Calcium Ion Concentration:} The second input $u_{n-1}=u(t_{n-1})$ captures the concentration of calcium ions \((\text{Ca}^{2+})\) at the time $t_{n-1}$, acknowledging its significant role in ion channel regulation. 
	\item \textbf{Rate of Change of Calcium Ion Concentration:} The third input $(u_n - u_{n-1})/\Delta t$, where $u_n = u(t_{n})$ and $\Delta t$ is the time step size, measures the rate of change in calcium ion concentration over time. This provides the network with more information about the dynamics of calcium levels. An increase or decrease in calcium concentration causes different reactions in the RyR channel.
\end{enumerate}
The output is designed as follows:
\begin{enumerate}
	\item \textbf{Rate of Change of The Open Probability of RyR Channel:} The output $(P_n - P_{n-1})/\Delta t$ is an approximation of the time derivative of $P(t)$ at time $t_{n-1}$. 
\end{enumerate}
Details of the fully connected Neural Network are summarized in the following table:
\begin{center}
\begin{tabular}{ |c|c|c|c| } 
\hline
Layers         & Shape & Parameters & Activation\\ \hline
Input layer    & 3     & -    & -\\ \hline
Hidden layer 1 & 200   & 800  & ReLU\\ \hline
Hidden layer 2 & 64    &12864 & ReLU\\ \hline
Hidden layer 3 & 16    &1040  & ReLU\\ \hline
Output layer   & 1     &17    & Linear\\ \hline
\end{tabular}
\end{center}
The neural network, with the number of neurons specified in the ``Shape" column, employs Rectified Linear Units (ReLU) as the activation function in its hidden layers. This choice is motivated by ReLU's ability to provide non-linear modeling capabilities and mitigate issues like the vanishing gradient problem. The numbers of weights and biases between different layers are listed in the ``Parameters" column. The output layer uses a linear activation function. The network is denoted as a scalar-valued function with three variables and parameters:
$$\mathcal{F}(\cdot,\cdot,\cdot, \Theta)$$
where $\Theta$ represents the undetermined weights and biases. 
A better discrete version of the DNN model is:
\begin{equation}
\frac{P_{n}-P_{n-1}}{\Delta t} = 
\mathcal{F}\left(P_{n-1}, u_{n-1},\frac{u_{n}-u_{n-1}}{\Delta t} , \Theta\right),
\end{equation}
where $n\geq 1$. The loss function is the mean square error:
\begin{equation}
L(\Theta) =  \frac{1}{N}
\sum\limits_{n=1}^N\left(
\mathcal{F}\left(P_{n-1}, u_{n-1},\frac{u_{n}-u_{n-1}}{\Delta t} , \Theta\right) -\frac{P_{n}-P_{n-1}}{\Delta t}
\right)^2.
\end{equation}
For training the network, we need to minimize the loss function by updating $\Theta$, which requires sufficient data. Once well-trained, to predict the open probability at the next time step, we can utilize the following formula:
\begin{equation}\label{P_n}
P_{n} = 
P_{n-1} + {\Delta t}\mathcal{F}\left(P_{n-1}, u_{n-1},\frac{u_{n-1}-u_{n-2}}{\Delta t} , \Theta\right).
\end{equation}
where $n\geq 2$.  For $n=1$, simply let  $u_{-1} = u_0$. Meanwhile, to ensure $P_{n}$ is between 0 and 1, values below 0 or above 1 are truncated.

\subsection{The Model Trained on Dataset from the ODE System}\label{tset_ode}
\begin{figure}[ht!]
	\centering
	\subfigure[]{%
		\includegraphics[width=0.4\linewidth]{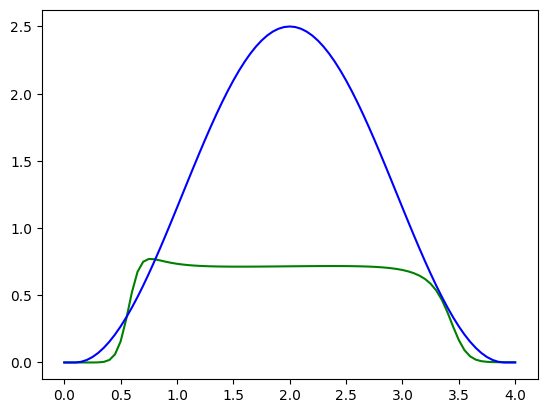}}%
	\hspace{4pt}%
	\subfigure[]{%
		\includegraphics[width=0.4\linewidth]{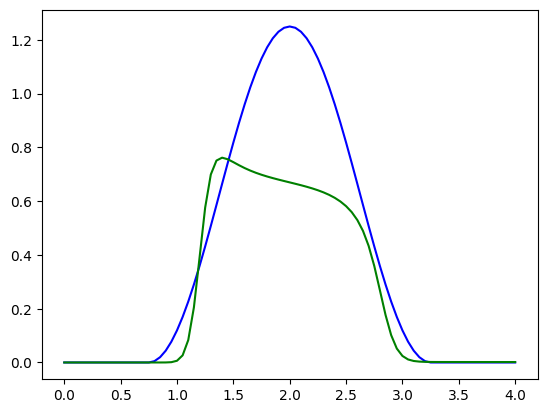}} 
    \subfigure[]{%
		\includegraphics[width=0.4\linewidth]{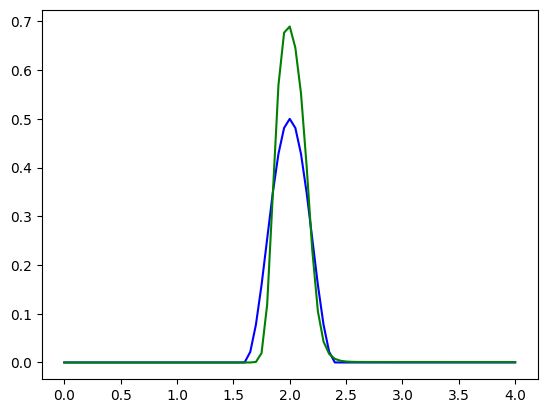}}%
	\hspace{4pt}%
	\subfigure[]{%
		\includegraphics[width=0.4\linewidth]{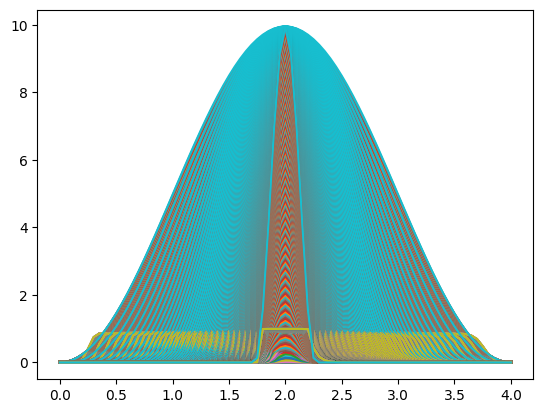}} 
  \caption[figures]{ Section \ref{tset_ode}. Training set from the ODE system. In all figures, the $x$-axis represents time, while the $y$-axis represents calcium concentration $u(t)$ and the open probability $P(t)$.
 From (a) to (c), the blue curves represent calcium signals, which fluctuate up and down over time with different magnitudes and durations; the green curves represent the corresponding open probability, ranging between 0 and 1. Figure (d) shows the collection of the entire training set, consisting of 26,000 pairs of $P(t)$ and $u(t)$. }
 \label{train1} 
\end{figure}

Training set is generated by the ODE system \eqref{ryr_ptu}-\eqref{ryr_ode} using Backward Euler method with time step size $\Delta t=0.05$. However, we only keep 
the values of open probability $P(t)$ and calcium concentration $u(t)$, $t\in [0,4]$. The parameters of equation \eqref{ryr_ode} are the same as in \cite{guan2022}. The initial values $(c_1, o, c_2)$ are set to $(1, 0, 0)$.
We choose 26,000 different calcium signals and define them as $\{u^k(t)\}_{k=1}^{26000}$.  The training set can be represented as
$$
\left\{\left(P_n^k, u_n^k, \frac{u_{n}^k-u_{n-1}^k}{\Delta t}\right),  \frac{P_{n}^k-P_{n-1}^k}{\Delta t}\right\}_{n=1,\cdots,80}^{k=1,\cdots,26000},
$$
where $P_n^k$ is the solution of equations \eqref{ryr_ptu}-\eqref{ryr_ode} with the $k$-th calcium signal $u^k(t)$ at time step $n$. $u^k(t)$ is a scaled and truncated cosine function centered at $t = 2$. Its magnitude ranges from 0.05 to 10, and its duration ranges from 0.5 to 4. The pairs of $(u^k(t),P^k(t))$ are shown in Figure \ref{train1}.

\begin{figure}[ht!]
\centering
	\subfigure[]{
		\includegraphics[width=0.4\textwidth]{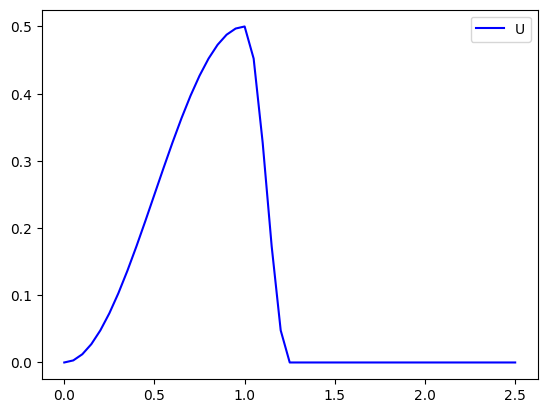}}
	\subfigure[]{
		\includegraphics[width=0.4\textwidth]{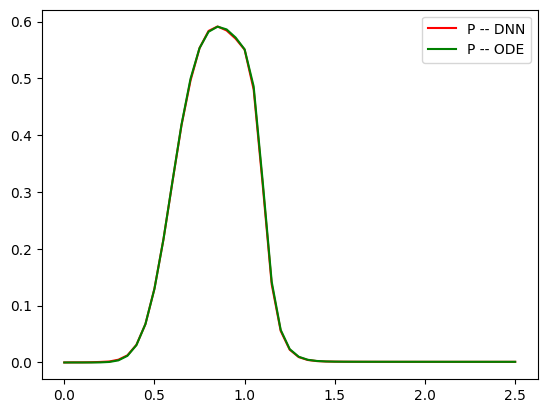}}
	\subfigure[]{ 
		\includegraphics[width=0.4\textwidth]{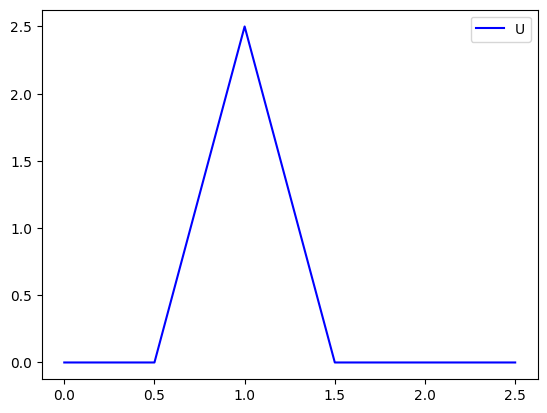}}
	\subfigure[]{
		\includegraphics[width=0.4\textwidth]{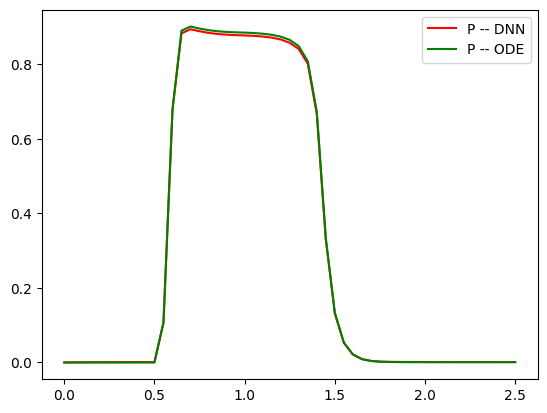}}
	\caption[figures]{ Section \ref{tset_ode}. Performance of the model on unseen data. (a) is an asymmetric calcium signal, and (b) is the response. (c) is a symmetric calcium signal, and (d) is its response. In both (b) and (d), the open probability $P$ generated by the neural network is shown in red, while the open probability from the ODE system is shown in green. The magnitudes of calcium signals are within the range of the training set.}
 \label{fig4}
\end{figure}

\begin{figure}[ht!]
\centering
	\subfigure[]{
		\includegraphics[width=0.4\textwidth]{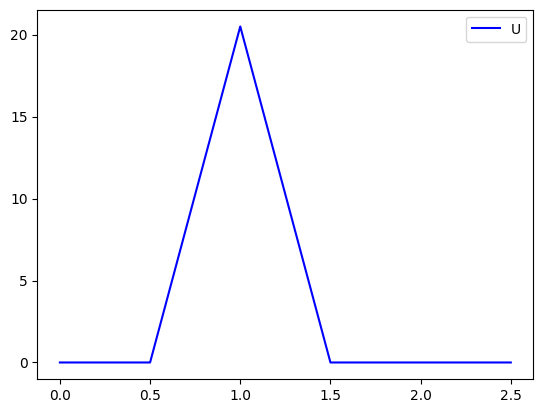}}
	\subfigure[]{
		\includegraphics[width=0.4\textwidth]{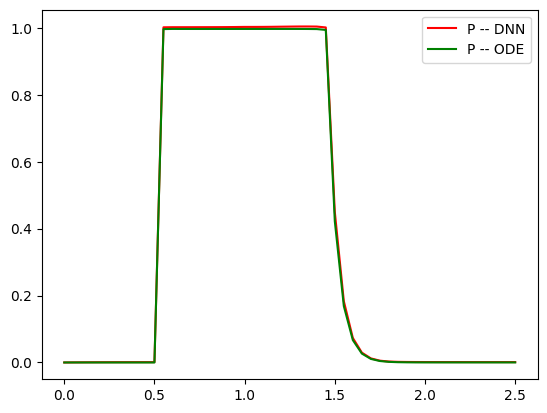}}
	\subfigure[]{ 
		\includegraphics[width=0.4\textwidth]{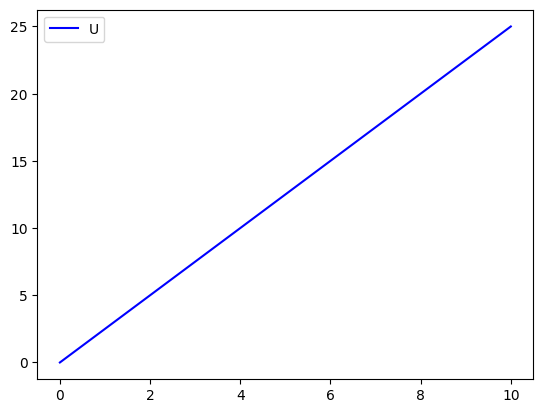}}
	\subfigure[]{
		\includegraphics[width=0.4\textwidth]{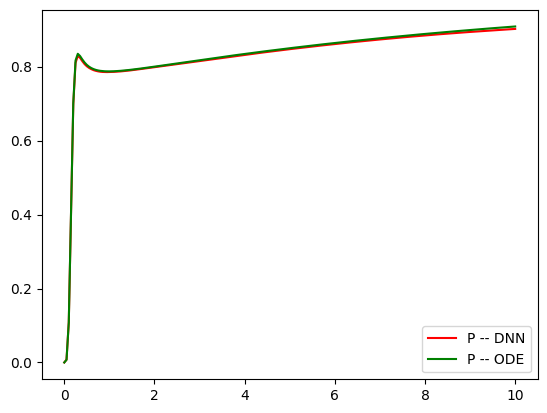}}
	\caption[figures]{ Section \ref{tset_ode}. Performance of the model on unseen data. (a) is a symmetric calcium signal, and (b) is the response. (c) is an asymmetric calcium signal, and (d) is its response. In both (b) and (d), the open probability $P$ generated by the neural network is shown in red, while the open probability from the ODE system is shown in green. The magnitudes of calcium signals are beyond the range of training set.}
 \label{fig5}
\end{figure}

To train the neural network, we set the batch size to 640, with 100 epochs, and 10\% validation data. We utilize the Adam optimizer, a method known for its efficiency and adaptive learning rate capabilities. After training, we test the model with unseen calcium signals, which are not in the training set. The results are given in Figures \ref{fig4}-\ref{fig5}. 

In Figure \ref{fig4}, the amplitudes of the calcium signals are 0.5 and 2.5. Both are within the range of the training set. The results match well with the open probability generated by the ODE system.

In Figure \ref{fig5}, the amplitudes of the calcium signals are 20 and 25, which are beyond the upper bound of the training set. The DNN model still works well. This proves the generalization ability of the DNN model.

However, in Figure \ref{fig5}(c)-(d), as calcium concentration increases, the open probability goes up then drops and increases again for both the ODE model and the DNN model. This is physiologically unreasonable. Study in \cite{bezprozvanny1991bell} shows that the calcium dependence of RyR channels is biphasic. A calcium-induced decrease occurs at high calcium concentrations. In \cite{gyorke1993ryanodine}, it is showed that the Ca$^{2+}$ sensitivity of RyR channels decreases during prolonged exposure to Ca$^{2+}$. In reality, the open probability should drop if the calcium concentration rises too high. The DNN model, trained by data from the ODE system, also inherits its drawbacks. In the next section, we will build a new DNN model from artificial data to show how to counter this disadvantage.

\subsection{The Model Trained on Artificial Dataset}\label{AT}
\begin{figure}[ht!]
	\centering
	\subfigure[]{%
		\includegraphics[width=0.45\linewidth]{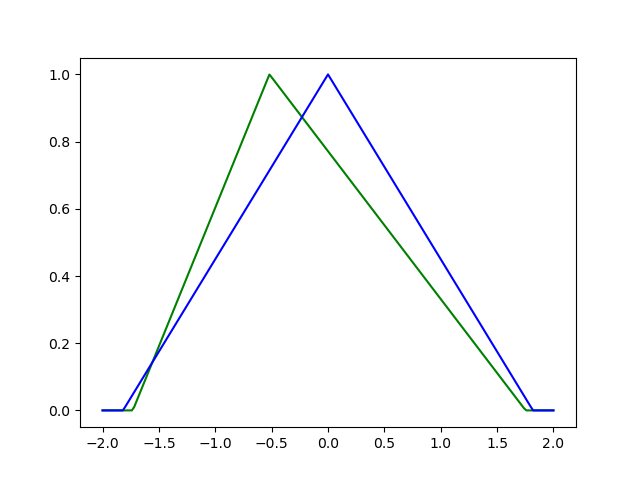}}%
	\subfigure[]{%
		\includegraphics[width=0.45\linewidth]{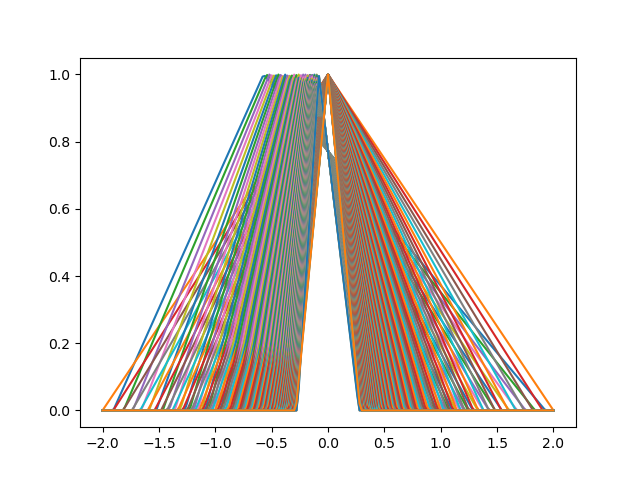}} 
	\caption[figures]{{\bf Training Set I} in Section \ref{AT}. The $x$-axis represents time, while the
$y$-axis represents calcium concentration $u(t)$ and the open probability $P(t)$. (a) shows one pair of $u(t)$ (in blue) and $P(t)$ (in green). (b) shows the whole training set, consisting of 121 pairs. }
 \label{train_hat_1}
\end{figure}

\begin{figure}[ht!]
	\centering
	\subfigure[]{%
		\includegraphics[width=0.45\linewidth]{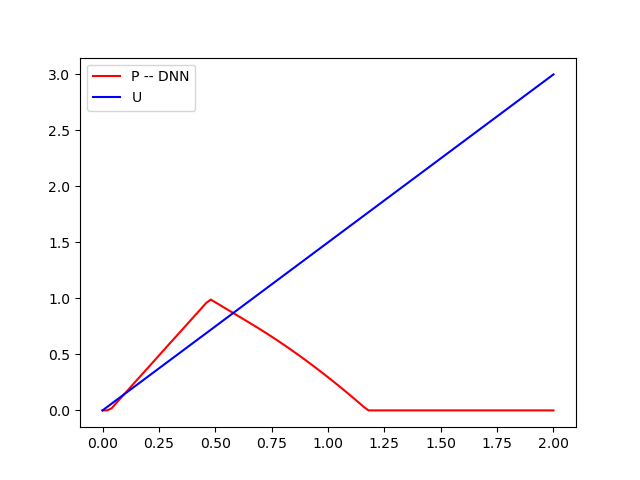}}%
	\hspace{0pt}%
	\subfigure[]{%
		\includegraphics[width=0.45\linewidth]{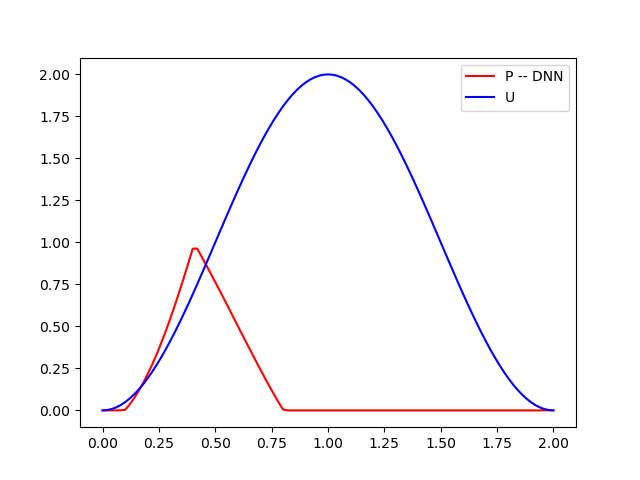}} 
	\caption[figures]{Model trained on {\bf Training Set I} in Section \ref{AT}. The $x$-axis represents time, while the $y$-axis represents calcium concentration $u(t)$ (in blue) and the open probability $P(t)$ (in red).}
  \label{train_hat}
\end{figure}

\begin{figure}[ht!]
		\centering
		\subfigure[]{%
			\includegraphics[width=0.45\linewidth]{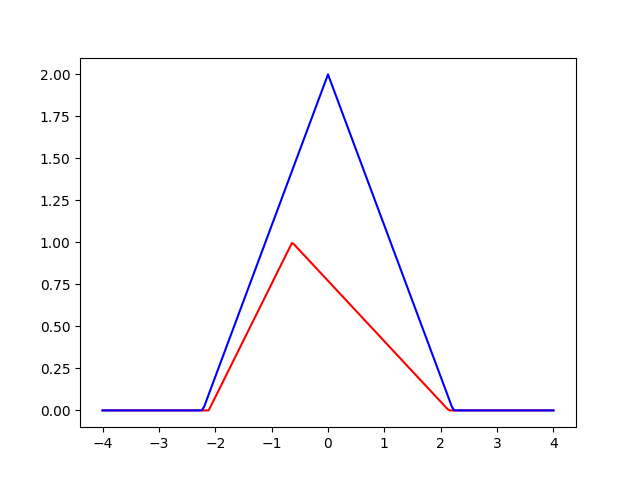}}%
		\hspace{4pt}%
		\subfigure[]{%
			\includegraphics[width=0.45\linewidth]{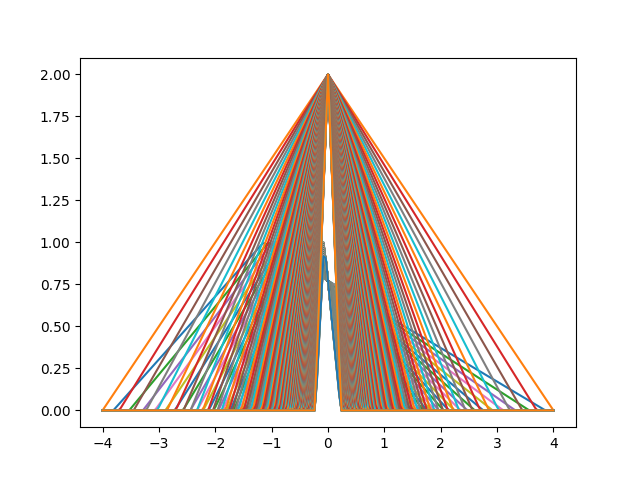}} 
			\caption[figures]{{\bf Training Set II} in Section \ref{AT}. The $x$-axis represents time, while the $y$-axis represents calcium concentration $u(t)$ and the open probability $P(t)$. (a) shows one pair of $u(t)$ (in blue) and $P(t)$ (in red). (b) shows the whole training set, consisting of 200 pairs.}
  \label{fig8}
\end{figure}

\begin{figure}[ht!]
		\centering
		\subfigure[]{%
			\includegraphics[width=0.45\linewidth]{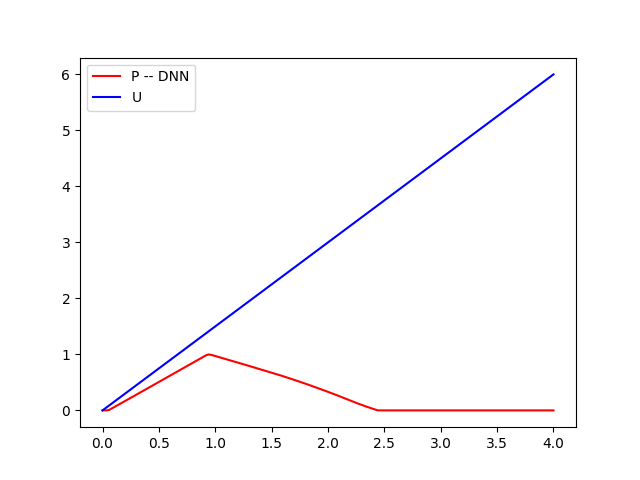}}%
		\hspace{4pt}%
		\subfigure[]{%
			\includegraphics[width=0.45\linewidth]{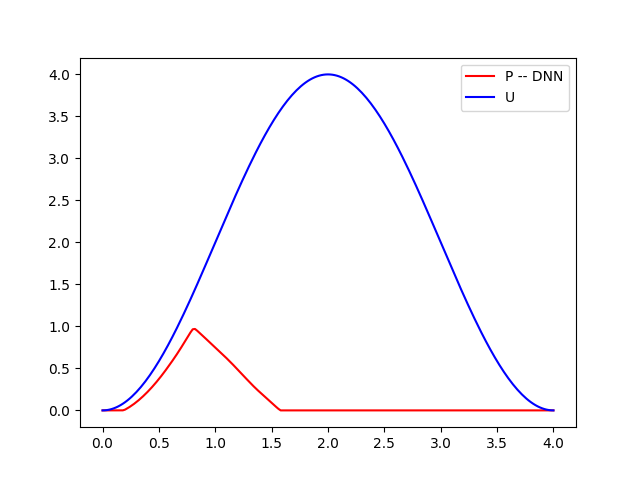}} 
   \caption[figures]{Model trained on {\bf Training Set II} in Section \ref{AT}. The $x$-axis represents time, while the $y$-axis represents calcium concentration $u(t)$ (in blue) and the open probability $P(t)$ (in red).}
   \label{fig9}
\end{figure}

The data in this section is created by following the patterns of the training set in Section \ref{tset_ode}. As in Figure 3(a)-(c), the open probability rises after the calcium concentration goes up, reaches its peak before the calcium signal does, and drops before the calcium signal drops.
However, we fix the amplitudes of the calcium signals. We construct two training sets and train the models on them. Details are summarized as follows.
\begin{enumerate}
	\item \textbf{Training Set I:} As in Figure \ref{train_hat_1}, there are 121 pairs of $u(t)$ and $P(t)$, where $t \in [-2, 2]$. The time step size is $\Delta t = 0.02$. The amplitudes of $u(t)$ and $P(t)$ are both 1.
    \item[] \textbf{Results:} From Figure~\ref{train_hat}(a), we see that unlike Figures~\ref{fig5}(c) and~(d), as calcium concentration increases, $P(t)$ initially goes up and then drops down. This means we can counter the drawback of the ODE system by modifying the training set. From Figure~\ref{train_hat}(b), the open probability turns down before $u(t)$ reaches 1, which is consistent with the training set.
	\item \textbf{Training Set II:} In Figure \ref{fig8}, there are 200 pairs of $u(t)$ and $P(t)$, where $t \in [-4, 4]$. The time step size is $\Delta t = 0.04$. The amplitude of $u(t)$ is 2. 
     \item[] \textbf{Results:} From Figure~\ref{fig9}(a), we see that the open probability $P(t)$ behaves similarly as in Figure \ref{train_hat}. However, it decreases when $u(t)$ is greater than 1, and the duration of $P(t)$ is longer.  
\end{enumerate} 
After training the neural network on different datasets, we obtain new models for the open probability of the RyR channel. By adjusting the training sets, we can influence the opening time of the channel and determine when the open probability decreases. These quantities are difficult to modify with the ODE model. In the next section, the DNN models are incorporated into the PDE system \eqref{diff_r}.  

\section{Numerical Method for the PDE-DNN Model}\label{sec4}
In this section, we propose a one-step implicit-explicit time-stepping scheme to solve the hybrid PDE-DNN model. This scheme is more stable compared with the one in \cite{guan2022}. In space, we use the piecewise linear finite element method to discretize the polar PDE system. However, the application of FEM to PDEs in polar coordinates is hindered by the singularity of the terms with coefficient $1/r$ when $r=0$. There are two ways to deal with the sigularity. One is multiplying $r$ on both sides of the PDE, this has been done successfully, an implementation can be found in the documentation of FreeFem++, see \cite{MR3043640}. However, we are going to the other direction. By treating the hypersingular integral using the Hadamard finite-part integral (see Guan et al., 2022), we can also recover the full accuracy of the FEM. Equation \eqref{diff_r} shows that the singularity only appears for
$u_e$, since $u_e$ is defined on $[0,L]$, $u$ and $b$ are defined on $[L,R]$. If we discretize equation \eqref{diff_r} only in space, and denote the resulting discretized solutions by $u_h$, $b_h$, and $u_{e_h}$, then multiplying by the test functions $v_h$ on $[L,R]$ and $w_h$ on $[0,L]$, by integration by parts and boundary conditions in \eqref{neu_r}, we obtain the following equations: 

\begin{equation}\label{dis1}
\left\{
\begin{aligned}
\left(\frac{\partial u_h}{\partial{t}}, v_h  \right) 
& +D_c
\left(\frac{\partial u_h}{\partial{r}},\frac{\partial v_h}{\partial{r}} \right) 
-D_c
\left(\frac{1}{r}\frac{\partial u_h}{\partial{r}}, v_h  \right) \\
&  = 
\left( f(u_h,b_h), v_h \right) 
+g_c(u_h)v_h|_{R} - g_e(u_h,u_{e_h})v_h|_{L},
 \text{ on } [L,R],\\
\left(\frac{\partial b_h}{\partial{t}}, v_h  \right) 
& +D_b
\left(\frac{\partial b_h}{\partial{r}},\frac{\partial v_h}{\partial{r}} \right) 
-D_b
\left(\frac{1}{r}\frac{\partial b_h}{\partial{r}}, v_h  \right) 
 = 
\left( f(u_h,b_h), v_h \right),
 \text{   on } [L,R],\\
\left(\frac{\partial u_{e_h}}{\partial{t}}, w_h  \right) 
&+D_c
\left(\frac{\partial u_{e_h}}{\partial{r}},\frac{\partial w_h}{\partial{r}} \right) 
-D_c
\left(\frac{1}{r}\frac{\partial u_{e_h}}{\partial{r}}, w_h  \right) \\
&= g_e(u_h,u_{e_h})w_h|_{L}, \text{  on } [0,L],
\end{aligned}
\right. 
\end{equation}
where $f(u_h,b_h) = K_b^-(b^0-b_h)-K_b^+ b_hu_h$, and $|_L$ or $|_R$ means fetching the values at $r=L$ or $r=R$. For simplicity, if we omit $v_h$ and $w_h$ in the inner products and right hand sides of equation 
\eqref{dis1}, then the implicit-explicit time stepping scheme is formulated as
\begin{equation}\label{dis2}
\left\{
\begin{aligned}
\frac{u_h^{n+1}-u_h^{n}}{\Delta{t}}
& +D_c
\frac{\partial u_h^{n+1}}{\partial{r}}
-D_c
\frac{1}{r}\frac{\partial u_h^{n+1}}{\partial{r}}\\
&  = 
f(u_h^{n+1},b_h^{n})
+g_c(u_h^{n+1})|_{R} - \bar{g}_e(u_h^{n+1},u_{e_h}^{n})|_{L},
 \text{ on } [L,R],\\
\frac{b_h^{n+1}-b_h^{n}}{\Delta{t}}
& +D_b
\frac{\partial b_h^{n+1}}{\partial{r}}
-D_b
\frac{1}{r}\frac{\partial b_h^{n+1}}{\partial{r}}
 = 
f(u_h^n,b_h^{n+1}),
 \text{ on } [L,R],\\
\frac{u_{e_h}^{n+1}-u_{e_h}^{n}}{\Delta{t}}
&+D_c
\frac{\partial u_{e_h}^{n+1}}{\partial{r}}
-D_c
\frac{1}{r}\frac{\partial u_{e_h}^{n+1}}{\partial{r}}= g_e(u_h^{n},u_{e_h}^{n+1})|_{L}, \text{ on } [0,L],
\end{aligned}
\right. 
\end{equation}
where $n\geq 0$, the term $\bar{g}_e$ is defined as
$$
\bar{g}_e(u_h^{n+1},u_{e_h}^{n})  := 
{C_{1}^e { P_n } (u_h^{n+1}-u_{e_h}^{n})} 
	+
{C_2^e\frac{u_h^{n+1}}{(K_s+u^{n+1}_h)u_{e_h}^n}} 
	-
{C_3^e(u_{e_h}^n-u_h^{n+1})}  
$$
and $P_n$ is the DNN model for RyR channel defined in \eqref{P_n}. Since the initial values are all constants, we let $u_h^0=u_h^{-1}$. The three schemes in equation \eqref{dis2} are separated, they can be calculated independently at step $n+1$, providing numerical solutions $u_h^n$, $u_{e,h}^n$ and $b_h^n$. 

Next, we need to deal with the singularity of the term:
$$
\left(\frac{1}{r}\frac{\partial u_{e_h}}{\partial{r}}, w_h  \right), \text{ where } r\in[0,L]. 
$$
The interval $[0,L]$ is divided into $M$ equal sized sub-intervals. Let $h=L/M$, and $r_i = (i-1)h,$ where $i=1,\cdots, M+1$. Then, we define $u_{e_h}$ as 
$$u_{e_h} = \sum\limits_{i=1}^{M+1}u_{e,i}\phi_i(r),$$
where $u_{e,i}$ is the numerical solution at each $r_i$ and $\phi_i(r)$ is the hat function with $\phi_i(r_i)=1$ for $1\leq i\leq M+1$. Each element of the matrix $A$ for this singular term is calculated as 
\begin{equation}\label{A_ij}
A_{i,j} = \left(\frac{1}{r}\phi_j'(r), \phi_i(r)  \right), \text{ where } i,j = 1,2,\cdots, M+1.
\end{equation}
For matrix $A$, the singularity occurs only if $i=1, j\leq 2$. To get $A_{1,1}$ and $A_{1,2}$, we employ the following definition of Hadamard finite part integral: 
\begin{equation}\label{Hint}
\int_{a}^{b} \frac{f(x)}{x - a} \, dx = \lim\limits_{\epsilon \to 0} \left( \int_{a+\epsilon}^{b} \frac{f(x)}{x - a} \, dx + f(a) \ln \epsilon \right),
\end{equation}
where $f(x)$ is a continuous function on $[a,b]$, see \cite{guan2022collocation}. From equation \eqref{A_ij} and the definition in equation \eqref{Hint}, we have
\begin{align*}
A_{1,1} &=  \int_0^h\frac{1}{r}\frac{-1}{h}\frac{h-r}{h} \, dr 
         = \frac{-1}{h^2}\int_0^h\frac{h-r}{r} \, dr  = \frac{1-\ln h}{h},\\
A_{1,2} &=  \int_0^h\frac{1}{r}\frac{1}{h}\frac{h-r}{h} \, dr 
         = \frac{1}{h^2}\int_0^h\frac{h-r}{r} \, dr  = \frac{\ln h - 1}{h}.
\end{align*}
All other elements of matrix $A$ are well defined and there are no singularities for the rest  parts of equation \eqref{dis2}. The optimal convergence rates are verified for steady state and time dependent diffusion equations of $u_e$ on $[0,L]$ in polar coordinates.

\section{Experiments for the PDE-DNN Models}
In this section, we will test the hybrid PDE-DNN model with three neural networks trained on different datasets. The PDE part is given in equation \eqref{ex_1} with boundary conditions in \eqref{ex_1b}. When solving the PDE system numerically by the scheme in equation \eqref{dis2}, replacing $P(t,u)$ by the neural network in equation \eqref{P_n}, we obtain the hybrid model. 
Parameters and initial values of the original PDE-ODE model are given in the following table and equations \eqref{ex_1}-\eqref{ex_1b}. This example shows the initiation, propagation, and disappearance of a calcium wave in a closed 2D cell or an infinitely long cylindrical axon.

\noindent{\bf Example 1}. This example is based on \cite{guan2022}, which generates Ca$^{2+}$ waves. Units are adjusted so that $t$ is in seconds, and $u$, $b$ and $u_e$ are in $\mu$M:
\begin{equation}\label{ex_1}
	\left\{
	\begin{aligned}
&{\partial_t u} -220\frac{\partial^2 u}{\partial{r}^2}-220\frac{1}{r}\frac{\partial u}{\partial{r}}   = f(b,u) 
		\quad \text{on } (1.5,\pi)\\
&{\partial_t b} -20\frac{\partial^2 b}{\partial{r}^2}-20\frac{1}{r}\frac{\partial b}{\partial{r}} = f(b,u) 
		\quad \text{on } (1.5,\pi)\\
& {\partial_t u_e} -220\frac{\partial^2 u_e}{\partial{r}^2}-220\frac{1}{r}\frac{\partial u_e}{\partial{r}} = 0  \quad \text{on }  (0,1.5)
	\end{aligned}
	\right. 
\end{equation}
where the reaction term is $f(b,u) = K_b^-(b^0-b)-K_b^+ bu$. Boundary condition for $b$ is $\partial_r b =0$, if $r=1.5$ or $r=\pi$. Other boundary conditions for $u$ and $u_e$ are: 
\begin{equation}\label{ex_1b}
\left\{
\begin{aligned}
D_c\partial_r u 
	&= C_3(1000-u)-\frac{C_2u}{1.8+u}-\frac{C_1u^2}{0.06^2+u^2}
	+g(r,t), \ \text{ if } r=\pi, \\
-D_c\partial_r u  
	&= 
	C_{1}^e P(t,u)(u_e-u)-C_2^e\frac{u}{(0.18+u)u_e}+C_3^e(u_e-u) 
	\ \text{ if } r=1.5,\\
D_c\partial_r u_e
	&= 
	C_{1}^e P(t,u)(u-u_e)+C_2^e\frac{u}{(0.18+u)u_e}-C_3^e(u_e-u) 
	\ \text{ if } r=1.5,
\end{aligned}
\right. 
\end{equation}
where $g(\pi,t)$ is the calcium influx or stimuli coming in at $r=\pi$:
$$
g(\pi,t) = 1200t^2(1-t)^2, \text{ if } t\in[0,1]; \quad g(\pi,t) = 0, \text{ if } t>1. 
$$
Suppose this model is for a 2D cell (a disc within a disc with the same center) as shown in Figure \ref{fig1}. The radii of the two circles are 1.5 and $\pi$. After transforming into polar coordinates and assuming there is no change in $\theta$, we obtain the 1D model in equation \eqref{ex_1}. The 1D model is defined on the interval $[0,1.5] \cup [1.5,\pi]$. The spatial mesh size for the interval $[0,1.5]$ is $1.5/40$, and the mesh size for the interval $[1.5,\pi]$ is $(\pi-1.5)/40$. We employ the numerical method developed in Section \ref{sec4} and obtain the following results.

\begin{figure}[ht!]
	\centering
	\subfigure[PDE-DNN Model, $\Delta\!t=0.036$]{%
		\includegraphics[width=0.4\linewidth]{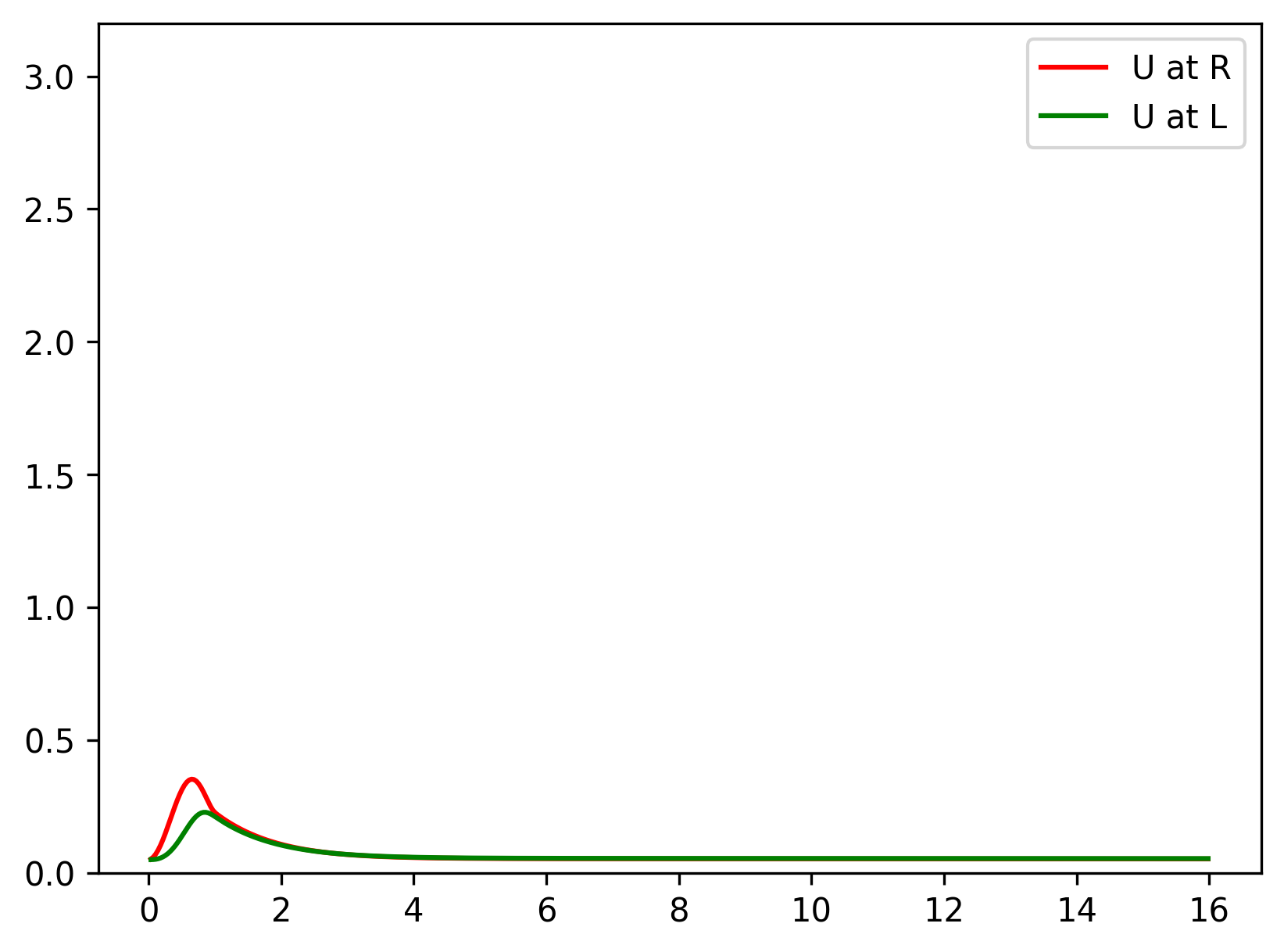}}%
	\hspace{4pt}%
	\subfigure[PDE-ODE Model, $\Delta\!t=0.036$]{%
	\includegraphics[width=0.4\linewidth]{png/no_RyR.png}} 
  \caption[figures]{Example 1. $P(t,u) = 0$. No calcium wave is generated. The calcium concentration $u$ at $r=R$ is in red, and $u$ at $r=L$ is in green. Calcium influx is injected at $r=R$, it starts at $t=0$ and ends at $t=1$. (a) is the model with DNN as RyR channel, (b) is the original ODE model. The results are the same. }
\label{fig10}
\end{figure}

\begin{table}
    \centering
    \begin{tabular}{cccccccc}
    \hline
    \multicolumn{8}{ c }{Parameters for the PDE-ODE Model} \\
    \hline
    $u(r,0)$ & $0.05$ && $b(r,0)$ & $37$ && $u_e(r,0)$ & $250$\\
    $c_1(0)$ & $0.994$ && $o(0)$ & $1.5721 \times 10^{-7}$ && $c_2(0)$ & $5.6625 \times 10^{-3}$ \\ 
    $b^0$ & $40$ && $K_b^-$ & $16.65$ && $K_b^+$ & $27$ \\ 
    $C_{1}^e$ & $0.829468$ && $C_2^e$ & $11000$ && $C_3^e$ & $0.038$ \\ 
    $C_1$ & $8.5$ && $C_2$ & $37.6$ && $C_3$ & $0.0045$ \\ 
    $k_a^{-}$ & $28.8$ && $k_b^{-}$ & $385.9$ && $k_c^{-}$ & 0.1 \\ 
    $k_a^{+}$ & 1500 && $k_b^{+}$ & 1500 && $k_c^{+}$ & 1.75 \\ 
    \hline
    \end{tabular}
\end{table}
To illustrate the importance of the RyR channel, we assume a scenario where $P(t,u) = 0$ constantly, which is equivalent to removing all RyR channels from the ER membrane. The result is shown in Figure \ref{fig10}.
Since RyR channel is responsible for releasing calcium ions from the ER into the cytoplasm in response to  the influx of calcium through plasma membrane. Without the RyR channel opening, the ER cannot release stored Ca$^{2+}$ into the cytoplasm. As a result, the generation of a calcium wave is prevented. The peak of the calcium concentration is below 0.5. Compared with the amplitudes of calcium waves in Figure \ref{fig11}, this is much lower.

\subsection{DNN Trained on the ODE Generated Data}
\begin{figure}[!ht] 
	\centering
 	\subfigure[PDE-DNN Model, $\Delta\!t=0.05$]{%
		\includegraphics[width=0.35\linewidth]{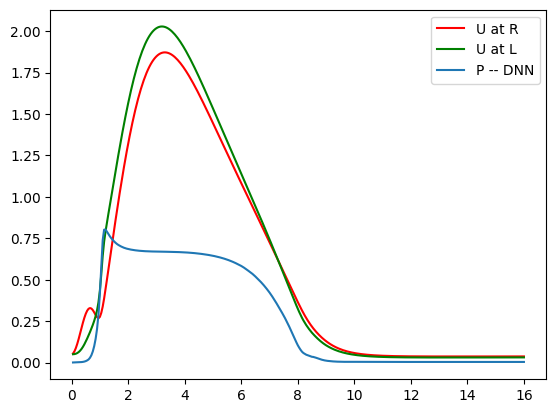}}%
	\hspace{4pt}%
	\subfigure[PDE-ODE Model, $\Delta\!t=0.05$]{%
		\includegraphics[width=0.35\linewidth]{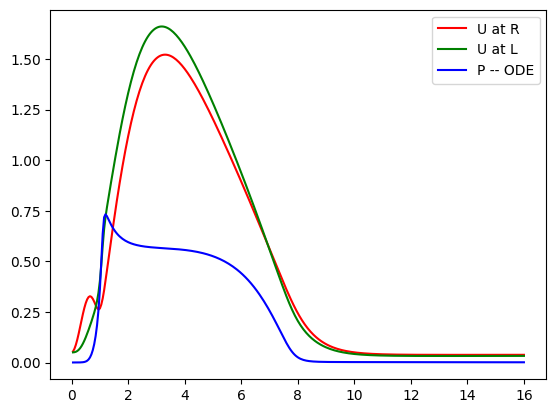}}
	\subfigure[PDE-DNN Model, $\Delta\!t=0.036$]{%
		\includegraphics[width=0.35\linewidth]{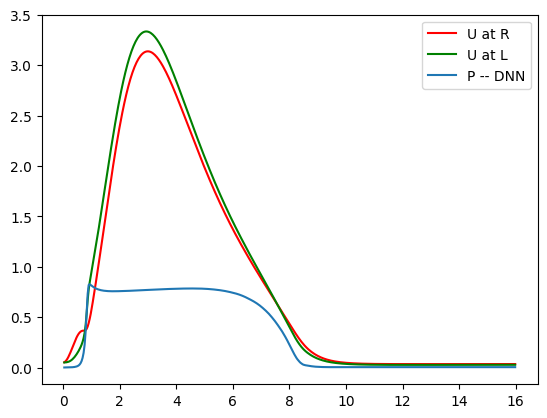}}%
	\hspace{4pt}%
	\subfigure[PDE-ODE Model, $\Delta\!t=0.036$]{%
		\includegraphics[width=0.35\linewidth]{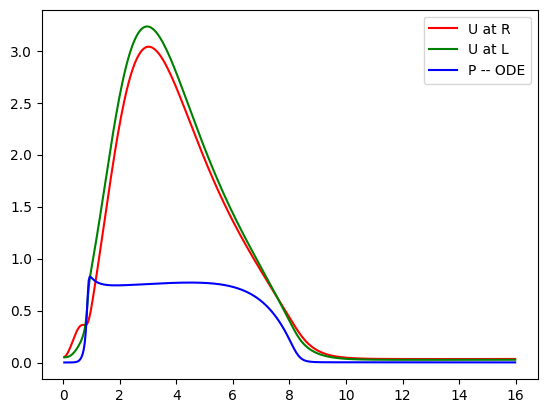}} 
	\subfigure[PDE-DNN, $\Delta\!t=0.025$]{%
		\includegraphics[width=0.35\linewidth]{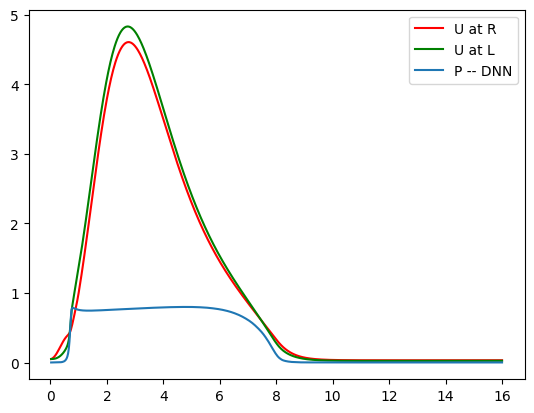}}%
	\hspace{4pt}%
	\subfigure[PDE-ODE, $\Delta\!t=0.025$]{%
		\includegraphics[width=0.35\linewidth]{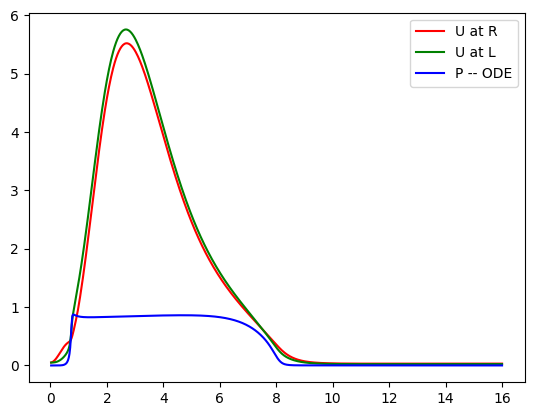}} 
	\subfigure[PDE-DNN, $\Delta\!t=0.0125$]{%
		\includegraphics[width=0.35\linewidth]{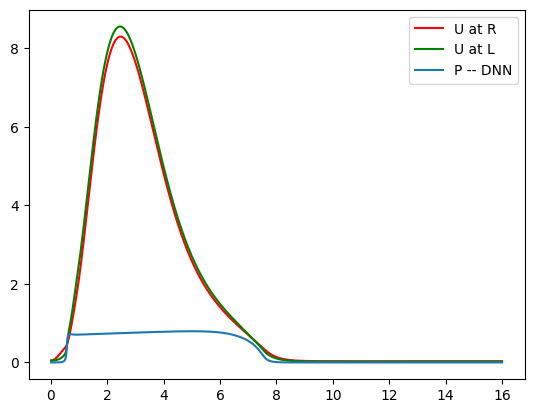}}%
	\hspace{4pt}%
	\subfigure[PDE-ODE, $\Delta\!t=0.0125$]{%
		\includegraphics[width=0.35\linewidth]{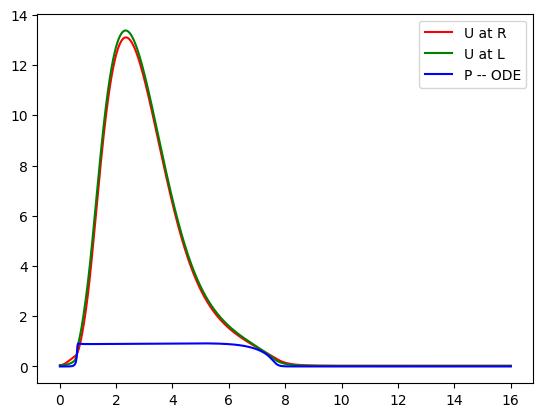}}
  \caption[figures]{Example 1. Section 5.1. The DNN model is trained on ODE-generated data. PDE-DNN models are listed on the left, including (a), (c), (e) and (g). PDE-ODE models are on the right side. The calcium concentration curve is green at $r=L$ (the ER membrane) and is red at $r=R$ (the plasma membrane). The open probability curve of the RyR channel is cyan for the DNN model and blue for the ODE model.}
  \label{fig11}
\end{figure}

If RyR channels function well, calcium waves can be generated.  As shown in Figure \ref{fig11}, firstly, the influx diffuses from $r=R$ (plasma membrane) to $r=L$ (ER membrane). Then, increasing Ca$^{2+}$ concentration at $r=L$ activates the RyR channel. Once the channel is open (which means the open probability $P$ increases), Ca$^{2+}$ ions are released from the ER. This causes the rise of calcium concentration in the cytosol (the region $[L,R]$), so the green and red curves go up beyond 0.5. Meanwhile, the open probability of the RyR channel must be dynamic, not static; otherwise, no wave will be produced.

\begin{figure}[ht!]\label{pde_ode}
	\centering
	\subfigure[$\Delta\!t=1/100$]{%
		\includegraphics[width=0.4\linewidth]{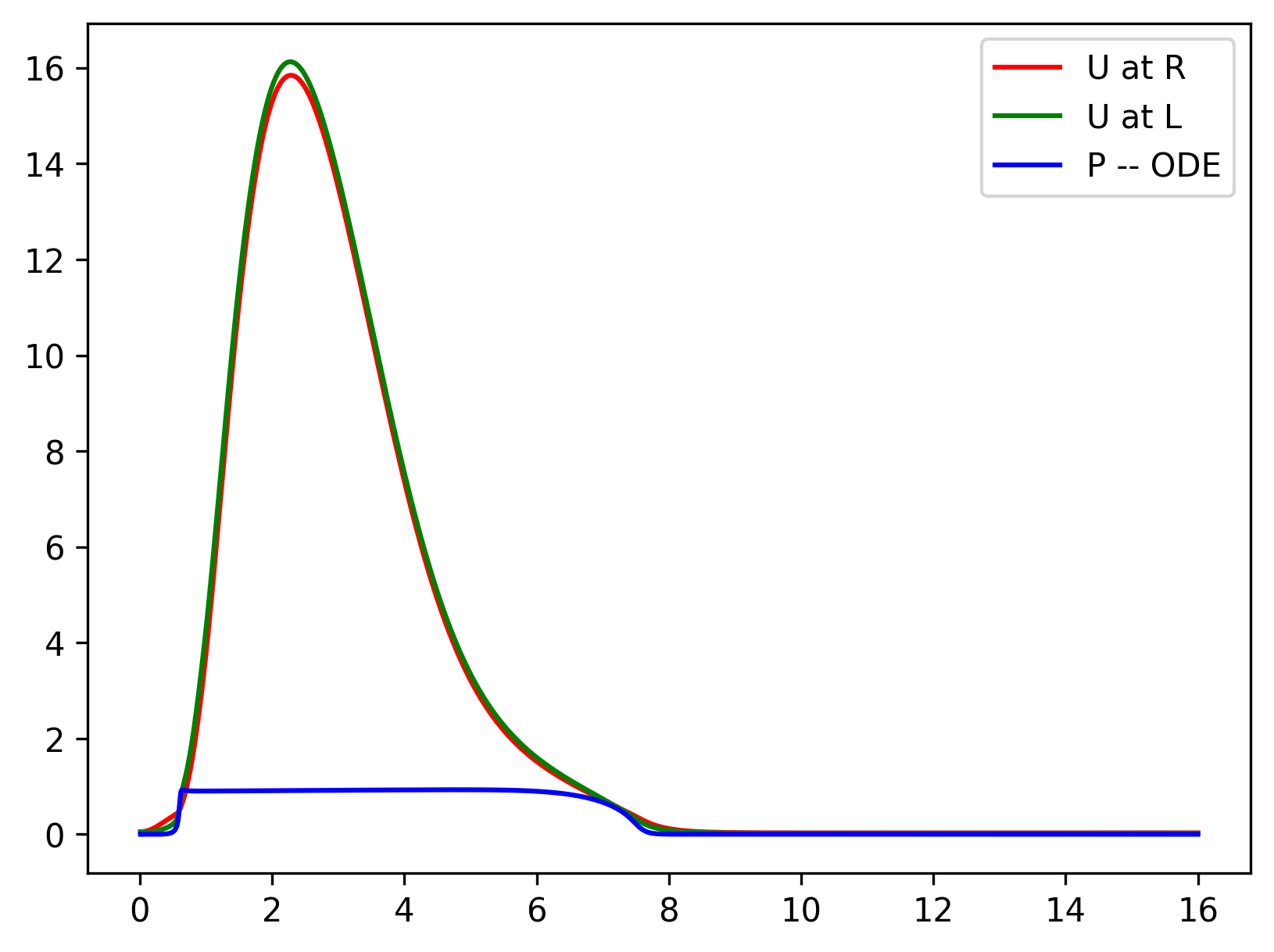}}%
	\hspace{4pt}%
	\subfigure[$\Delta\!t=1/500$]{%
		\includegraphics[width=0.4\linewidth]{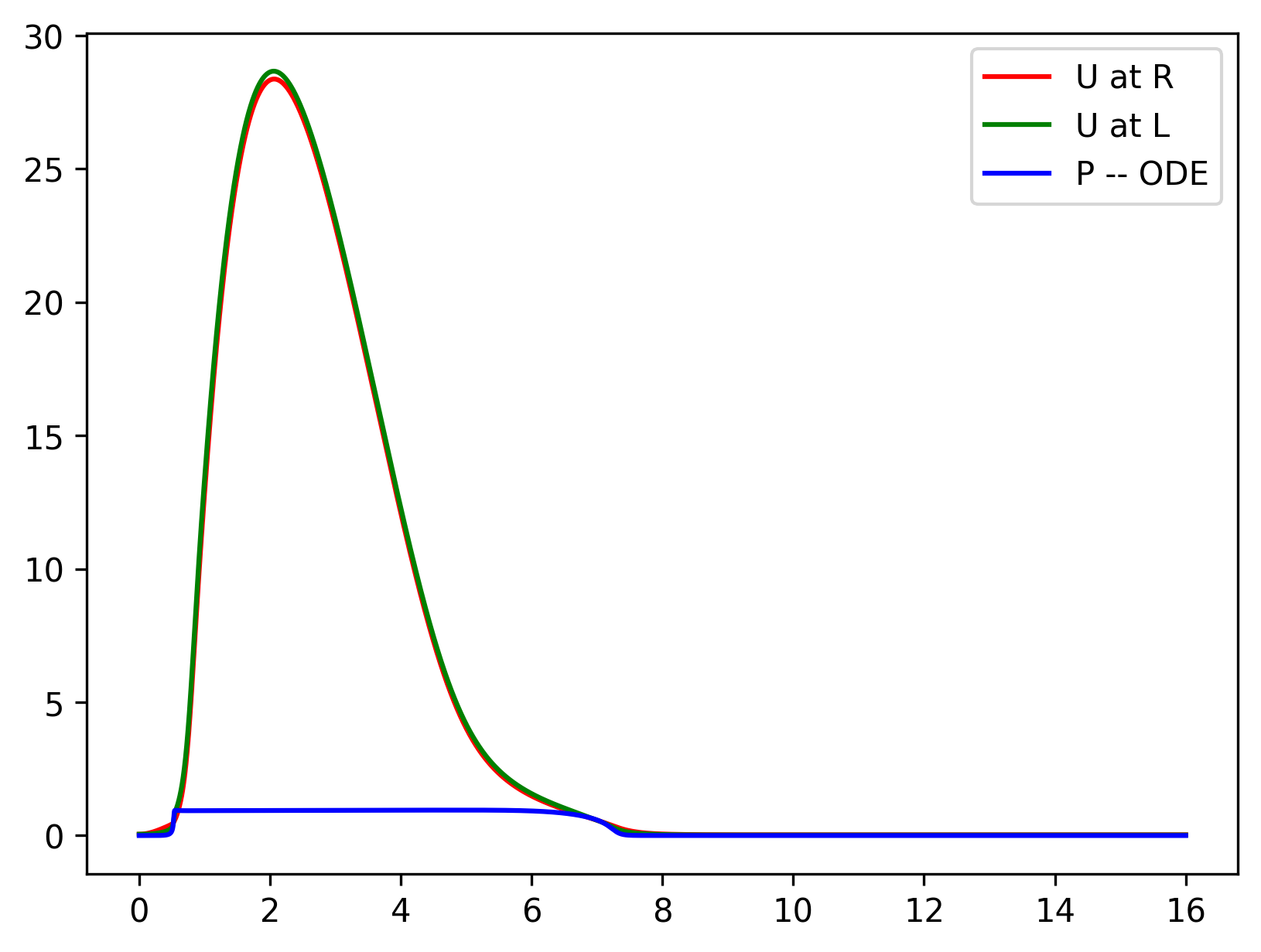}} 
	\subfigure[$\Delta\!t=1/2500$]{%
		\includegraphics[width=0.4\linewidth]{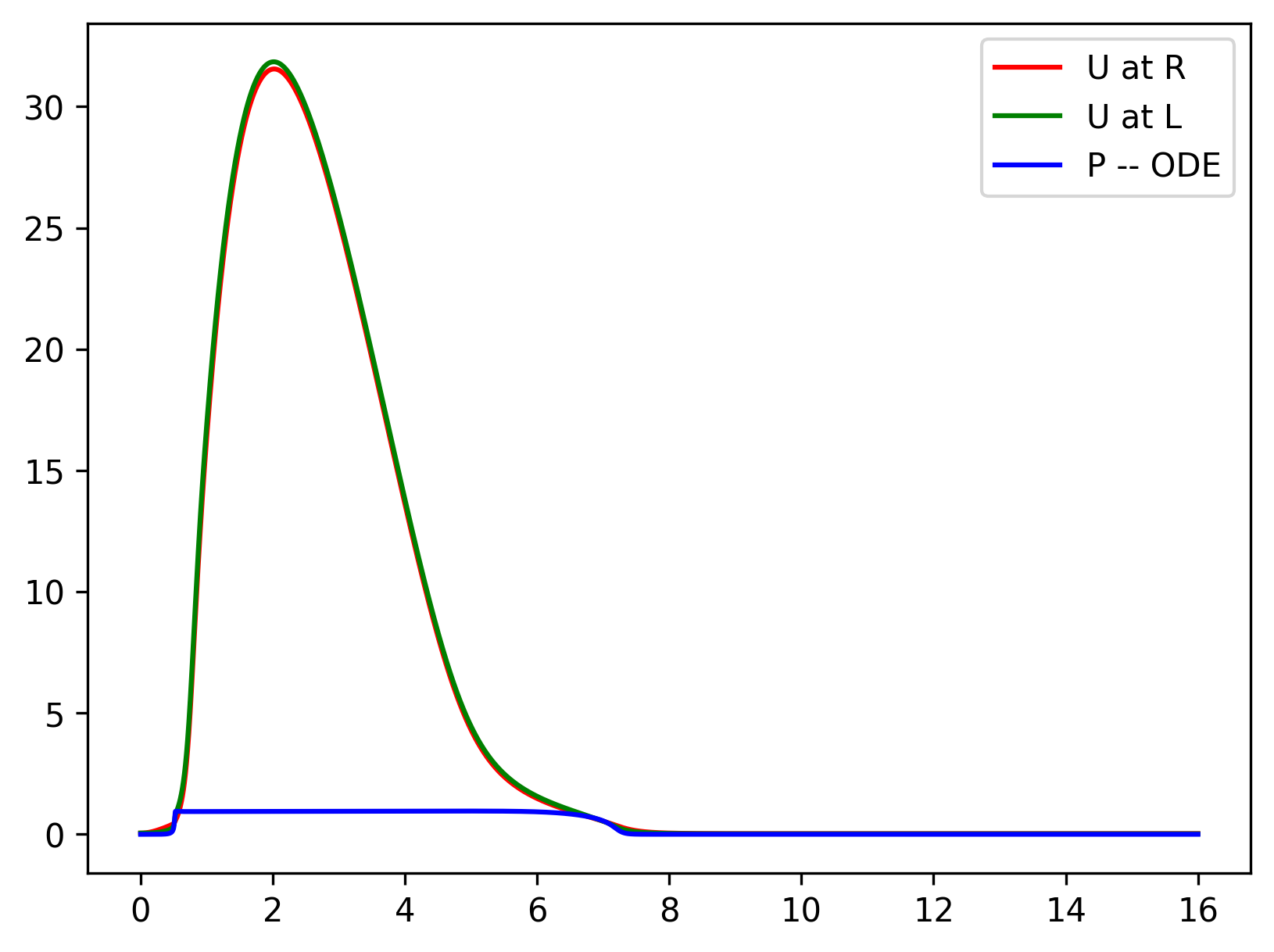}}%
	\hspace{4pt}%
	\subfigure[$\Delta\!t=1/12500$]{%
		\includegraphics[width=0.4\linewidth]{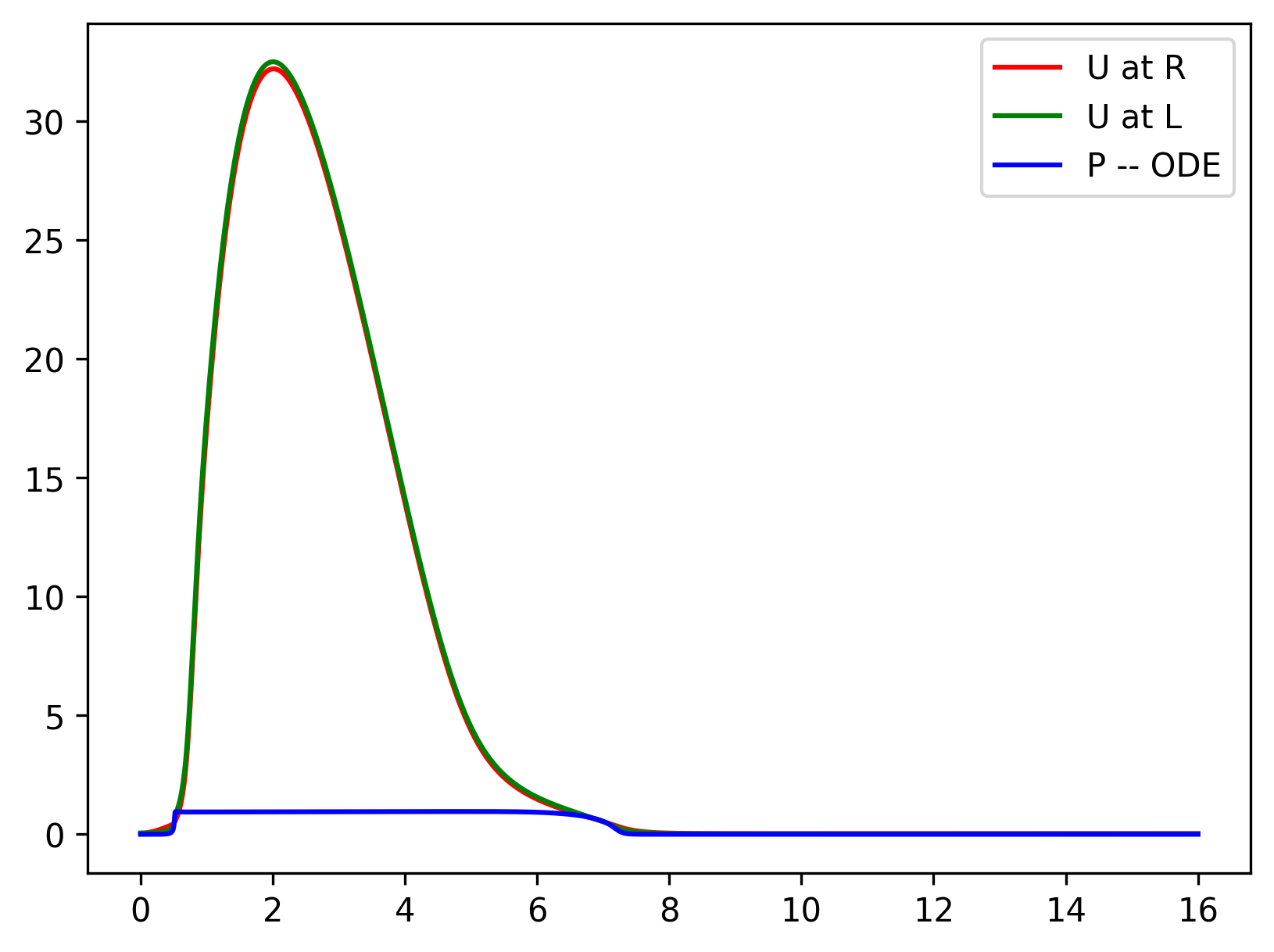}}
  \caption[figures]{Example 1. Section 5.1.  How can the PDE-ODE model achieve good accuracy? All figures here are for the PDE-ODE model. The amplitude of the calcium wave stabilizes around 32.}
  \label{fig12}
\end{figure}

By comparing the numerical results in Figure \ref{fig11}, we see that when $\Delta t$ is 0.05 or 0.036, the amplitudes from both models are close. This is because the DNN model is trained using data with a time step size of 0.05, so it performs well around 0.05 only. When $\Delta t$ is as small as 0.0125, the two models still capture a similar duration of the calcium wave, but the amplitudes are far apart. This is expected. However, if we compare the solutions in Figure \ref{fig11} with the accurate solution in Figure \ref{fig12}(d), there is still a significant gap. If $\Delta t$ is not small enough, neither the PDE-DNN model nor the PDE-ODE model can be accurate enough. The amplitude of the calcium wave from the PDE-ODE model stabilizes around 32 with a very small time step size $\Delta t=1/2500$, see Figure \ref{fig12}.
This disadvantage mainly stems from the ODE system in equation \eqref{ryr_ode}. It is not suitable to be used in a PDE model.

In the real world, we care more about the amplitude and duration of the calcium wave, not the shape of the wave. It's obvious that the PDE-ODE model is not a good one to use for simulation. It is too sensitive to time step size. Also, the ODE model is not physiologically reasonable as shown in Figure \ref{fig5}(c)-(d). By reducing the order of the ODE system, we show that DNN has the ability to learn the dynamics of the RyR channel. So, if there is reasonable data, we can get the DNN model directly from the data instead of building an ODE system. Next, we will combine the DNN models, which are trained on artificial datasets, with the PDE to run the simulation.

\subsection{DNN Trained on the Artificial Datasets}
There is no obvious way to improve the ODE model or adjust the amplitude or duration of the calcium wave from the PDE-ODE model. Maybe it's time to consider using new open probability models based on DNN. Since, it's too costly to use the PDE-ODE model. To stabilize the amplitude of the calcium wave from the PDE-ODE model, we need to use a very small time step. However, the PDE-DNN model is much less sensitive to time step size and saves a lot of runtime in total, if we don't count the time spent for training the DNN model. 
Actually, the artificial training sets are very small compared to the ODE-generated training set. All DNN models mentioned below are trained on artificial datasets.

To clearly show the calcium wave, we reduce the influx to 
$$
g(\pi,t) = 600t^2(1-t)^2, \text{ if } t\in[0,1]; \quad g(\pi,t) = 0, \text{ if } t>1. 
$$
By doing so, the calcium concentration at the ER membrane will be higher for the DNN model trained on 
{\bf Training Set I}. If not, the calcium concentration at $r=R$ will be higher, since the stronger influx vanishes slower. Of course, the influx can't be too small; otherwise, the RyR channel can't be activated.

\begin{figure}[ht!]
	\centering
	\subfigure[$\Delta\!t=1/80$]{%
		\includegraphics[width=0.4\linewidth]{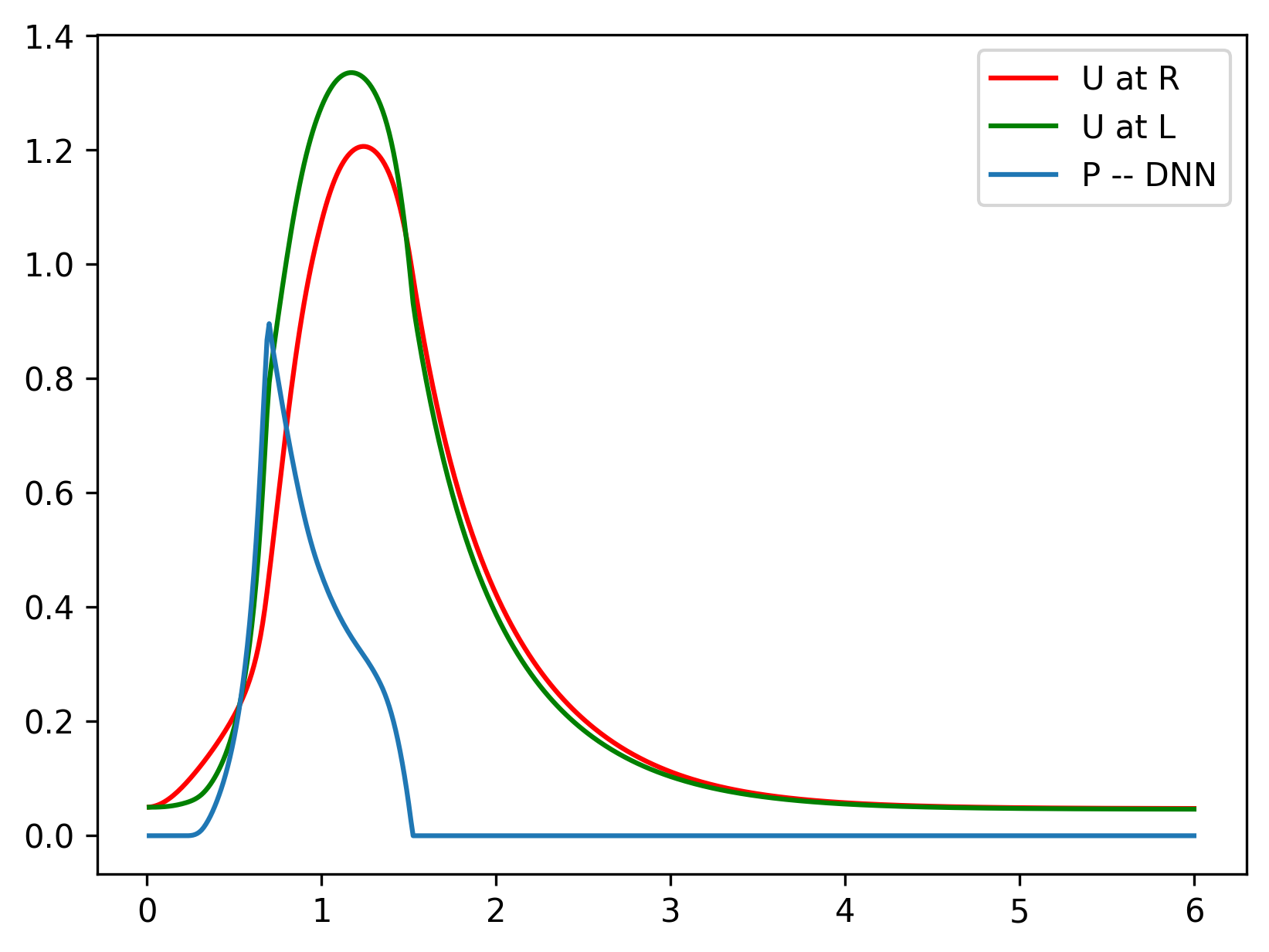}}%
	\hspace{4pt}%
	\subfigure[$\Delta\!t=1/160$]{%
		\includegraphics[width=0.4\linewidth]{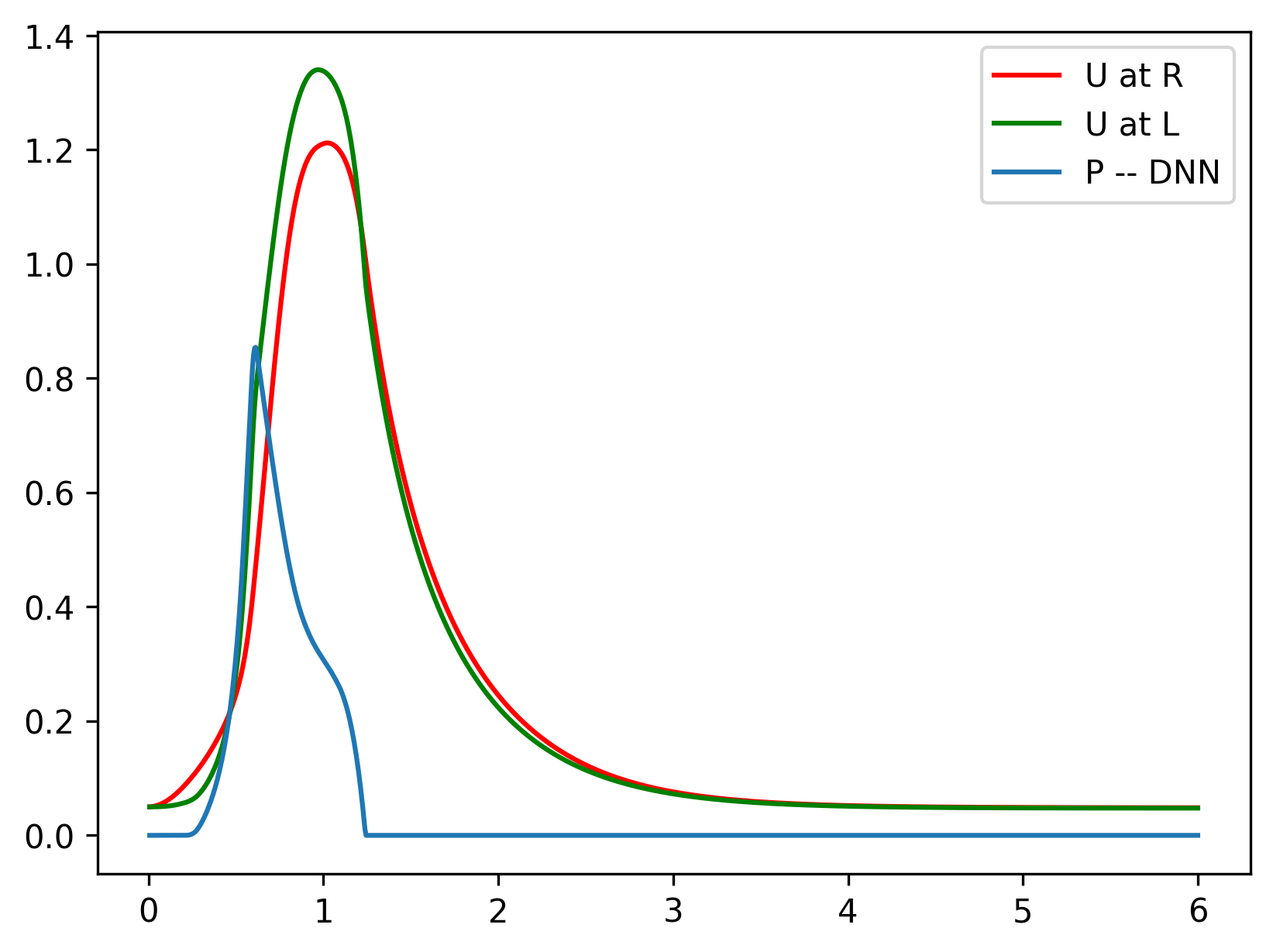}} 
  	\caption[figures]{Example 1. Section 5.2. Simulation of the PDE-DNN model with different time step sizes, when the DNN model is trained on {\bf Training Set I} in Section 3.3. }
   \label{f13}
\end{figure}

Figure \ref{f13} shows the numerical results for the PDE-DNN model, where the DNN is trained on {\bf Training Set I} in Section 3.3. For $\Delta t = 1/80$ and $\Delta t = 1/160$, the amplitudes are very close. This is the major difference compared to the models in Section 5.1. The durations of the two simulations are also close, though Figure \ref{f13}(b) is ``thinner". Therefore, we are fine using $\Delta t = 1/80$ to run the simulation.

Figure \ref{f13} shows similar results for the PDE-DNN model, where the DNN is trained on {\bf Training Set II} in Section 3.3. The new model is much less sensitive to the time step size. This property is particularly desired for large-scale problems.

Our tests show that the amplitudes of the calcium waves from the PDE-DNN model are proportional to the magnitudes of the calcium signals in the training sets. Similarly, the duration of the calcium waves is proportional to the average duration of the calcium signals in the training sets. Therefore, we can adjust the amplitudes and duration of the calcium waves by modifying the training set.

\begin{figure}[ht!]\label{train_hat_2_2}
		\centering
		\subfigure[$\Delta\!t=1/80$]{%
			\includegraphics[width=0.4\linewidth]{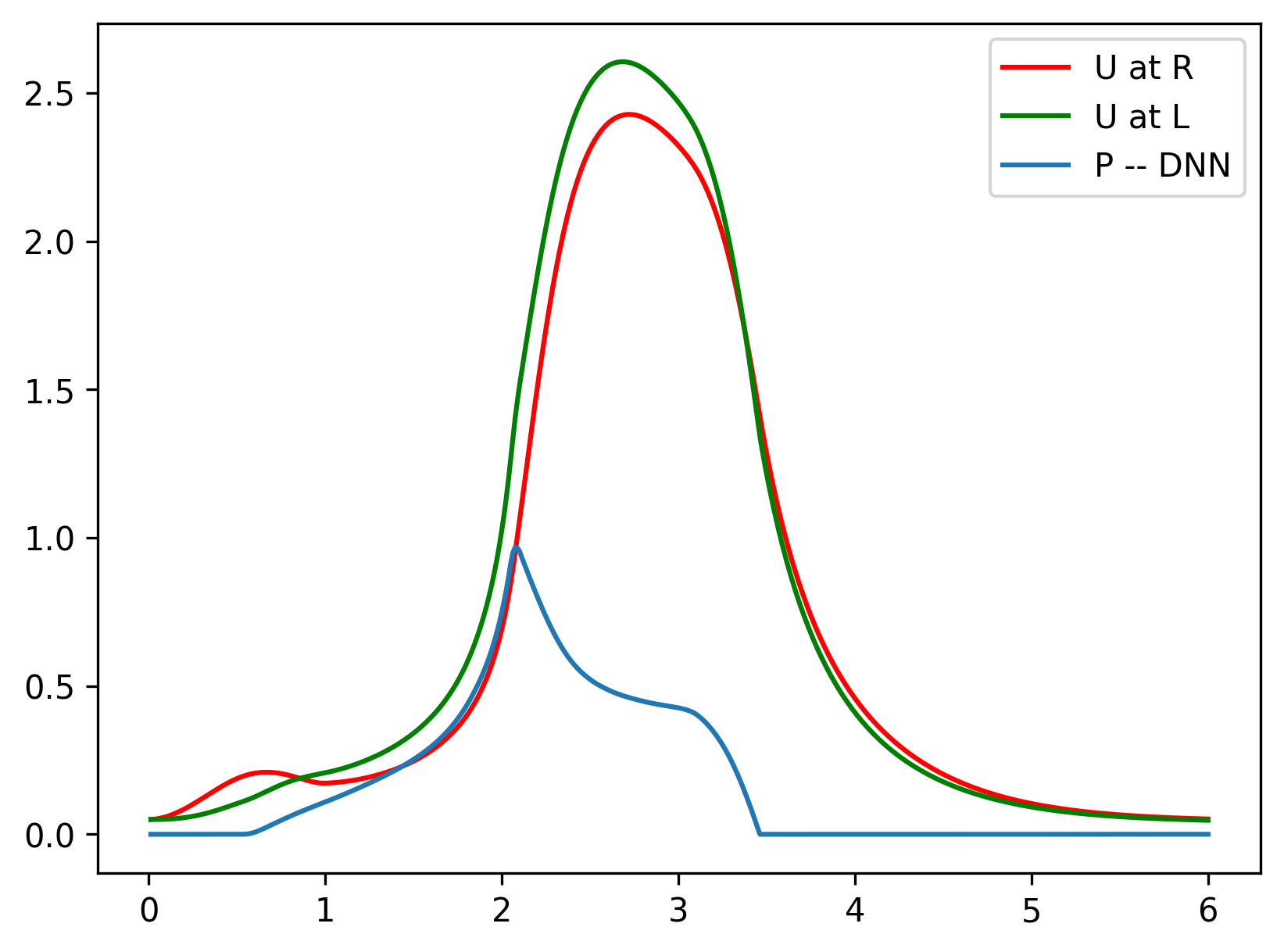}}%
		\hspace{4pt}%
		\subfigure[$\Delta\!t=1/160$]{%
			\includegraphics[width=0.4\linewidth]{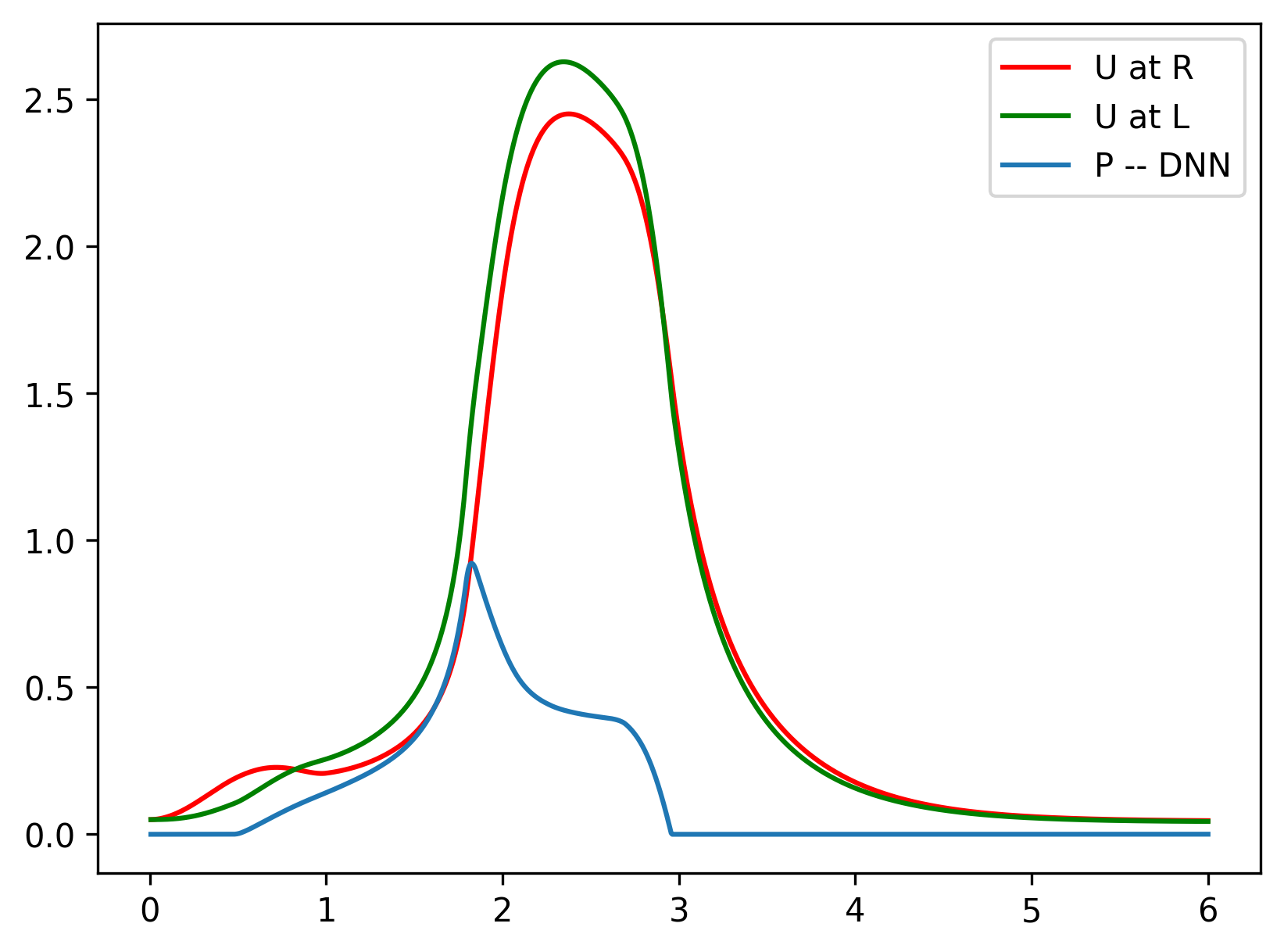}} 
            \caption[figures]{Example 1. Section 5.2. Simulation of the PDE-DNN model with different time step sizes, when the DNN model is trained on {\bf Training Set II} in Section 3.3. }
            \label{f14}
\end{figure}

\section{Conclusion}
In this paper, we build a hybrid PDE-DNN model for calcium dynamics in neurons. The PDE is simplified to a 1D diffusion-reaction system in polar coordinates, coupled via an interface point. We then design a stable time-stepping scheme and employ the Hadamard finite part integral to counter the singularity in the polar PDE. In this 1D calcium dynamic system, the ER membrane and plasma membrane are reduced to two boundary points, which makes it much easier and faster to prototype different models of ion channels.

The ODE system in the original PDE-ODE model is reduced to one equation with a single unknown by using a fully connected DNN. This technique can also be applied to other Markov chain models if the final output (like $P$ in equation \eqref{ryr_ptu}) of the ODE system is fewer or just one value, i.e., reducing the order of the original system.

The DNN model shows correct predictions but also inherits the drawbacks of the ODE system. To overcome the problems caused by the ODE system, we propose training the DNN on artificial but reasonable datasets, which lead to a better and adjustable model for RyR channels. The PDE-DNN model shows superior properties, such as insensitivity to time step size and the ability to easily change the amplitude and duration of the calcium waves. The DNN structure can also be used to model the open probability of other types of ion channels, such as Voltage-Gated Ion Channels, given experimental or reasonable artificial datasets.


\bibliographystyle{plain}

\bibliography{References}

\end{document}